\documentclass[11pt]{amsart}
 
\usepackage{graphicx}
\usepackage{amssymb,epsf}
\usepackage{amscd}\usepackage[dvips]{epsfig} 

\usepackage{amssymb} 
\usepackage{epsf} 
\usepackage{amsmath,epsf} 
\usepackage[latin1]{inputenc} 
\usepackage{amsfonts} 
\usepackage[dvips]{epsfig} 
\def\Young#1{\vbox{\smallskip\offinterlineskip 
    \halign{&\vbox{##}\kern-\Thickness\cr #1}}} 
 
\newdimen\Squaresize \Squaresize=12pt 
\newdimen\Thickness \Thickness=.3pt 
\newdimen\Correction \Correction=7pt 
 
\def\Vide#1{\hbox{ 
       \vbox to \Squaresize{\vss 
          \hbox to \Squaresize{\hss#1 \hss}\vss} 
    \hskip-\Correction} 
   \kern-\Thickness}

\def\Carre#1{\hbox{\vrule width \Thickness 
   \vbox to \Squaresize{\hrule height \Thickness\vss 
      \hbox to \Squaresize{\hss#1\hss} 
   \vss\hrule height\Thickness} 
   \unskip\vrule width \Thickness} 
   \kern-\Thickness}

\def\Box#1{\Carre{$\scriptstyle#1$}}

\newtheorem{theorem}{Theorem}[section]
\newtheorem{corollary}[theorem]{Corollary}
\newtheorem{prop}[theorem]{Proposition}
\newtheorem{lemma}[theorem]{Lemma}
\theoremstyle{definition}
\newtheorem{defn}{Definition}
\newtheorem{remark}{Remark}
\newtheorem{exem}[theorem]{Example}

\def\NCSI{{NCSym}}
\def\NSaNC{{\mathcal I}}
\def\H{{\mathcal H}}
\def\Sym{{Sym}}
\def\bbQ{{\mathbb Q}}
\def\la{{\lambda}}
\def\NSym{{NSym}}

\def\coeff{{\big|}}
\def\h{\hbox{{\bf h}}}
\def\e{\hbox{{\bf e}}}
\def\M{\hbox{{\bf m}}}
\def\cideal{{I}}
\def\R{{\bf R}}
\def\Sym{{\Lambda}}
\def\L{\mathcal{L}}
\def\C{{\mathbb{Q}}}
\def\A{\mathcal{A}}
\def\Sym{{Sym}}
\def\NCHar{{NCHar_n}}
\def\MH{{MHar}}
\def\haus{{\partial}}
\def\twist{{d}}
\def\twistCI{{\polys / \NCSIdeal}}
\def\s{\mathfrak{S}}
\def\polys{{\C\langle X_n \rangle}}
\def\NCSIdeal{{\langle \NCSI^+_n \rangle}}

\def\SymIdeal{{\langle \Sym^+_n \rangle}}
\def\vbar{\sep \, }
\def\dop{{(\![}}
\def\dcp{{]\!)}}
\def\shuf{{\sqcup\!\sqcup}}
\def\opsi{{\overline{\psi}}}
\def\sep{{,}}
\def\attach{{|}}

\title[Noncommutative Invariants and Coinvariants of $\s_n$]{Invariants and Coinvariants of the Symmetric Group in Noncommuting Variables}
\author[N.~Bergeron]{Nantel Bergeron}
\address{Mathematics and
Statistics, York University, 4700 Keele St., Toronto, ON, M5A 4T5, Canada}
\email{bergeron@mathstat.yorku.ca}
\email{rosas@mathstat.yorku.ca}
\email{zabrocki@mathstat.yorku.ca}
\author[C.~Reutenauer]{Christophe Reutenauer}
\address{LaCIM,
Un\-iver\-sit\'e du Qu\'e\-bec \`a Mont\-r\'eal,
Case Post\-ale 8888, suc\-cur\-sale Centre-ville,
Mont\-r\'eal (Qu\'e\-bec) H3C 3P8, Canada}
\email{christo@math.uqam.ca}
\author[M.~Rosas]{Mercedes Rosas}
\author[M.~Zabrocki]{Mike Zabrocki}

\begin{document}
\begin{abstract}
We introduce a natural Hopf algebra structure on the space of noncommutative
symmetric functions which was recently studied as a vector space by
Rosas and Sagan \cite{RS}. The bases for this algebra are indexed
by set partitions. We show that there exist a natural inclusion of the Hopf
algebra of noncommutative symmetric functions studied in \cite{thibon}
in this larger space.  We also consider this algebra as a subspace of
noncommutative polynomials and use it to 
understand the structure of the spaces of harmonics and coinvariants
with respect to this collection of noncommutative polynomials.
\end{abstract}

\maketitle

\begin{section}{Introduction}
In the commutative world  there are two
constructions of the Hopf algebra of symmetric functions;
the more classical one as the invariants of the symmetric group
on the polynomial ring, the other is the commutative free algebra
generated by one element in each degree (e.g. \cite{Mac} vs. \cite{S}).
These two constructions lead to the same algebra.

In the noncommutative world, these perspectives
lead to two very different
algebras; a free algebra with one generator at each degree, $\NSym$
(see for instance \cite{thibon} and
references therein), and the algebra of noncommutative invariant
polynomials, $\NCSI$ (studied in \cite{BeCo}, \cite{RS}, \cite{W}).  
These algebras are clearly not isomorphic
since the index set of a homogeneous basis of
$\NSym$ is the set of compositions, with dimension
$2^{n-1}$ at degree $n$, and  $\NCSI$ is indexed by the
set partitions, with dimension given by the Bell numbers for
each graded component.

An obvious question that was first posed in the work of \cite{RS}
was to understand the connection between the two algebras.
In this paper we present an incredibly beautiful relationship
between them.
We answer this question by first introducing a natural
Hopf algebra structure on the space of noncommutative
polynomial invariants.  The structure imposed by having both a product
and a coproduct is much richer than just the algebra structure
alone.  This places so many conditions on a Hopf algebra embedding
from $\NSym$ to $\NCSI$
that there is only one solution given a set of generators of $\Sym$.

In the development of the proof that the embedding is injective,
we compute the determinant of a combinatorial matrix indexed
by compositions.  It arises that the determinant is expressed as a
product of the number of permutations with
no global descents.  This is a surprising fact since these
numbers also happen to be the number of free generators/primitives
of the Malvenuto-Reutenauer Hopf algebra of permutations (see \cite{AS}, \cite{Co}).

In section \ref{sec:algebra} we digress and mention several
other relationships between $\NCSI$ and $\NSym$
by examining some quotients and
embeddings of the algebra structure.  In particular,
we show that, as a graded algebra, $\NCSI$ in two variables
is isomorphic to $\NSym$.

An interesting question that is natural to consider once
one understands $\NCSI$ as the space of invariants
is to try to understand the `coinvariants.'
An important classical theorem of Chevalley \cite{C}, and later 
extended to other finite reflection groups by Steinberg \cite{steinberg},
says that the ring of polynomials is isomorphic as an
$\s_n$-module to the tensor product of its invariants
times its coinvariants.
We next ask ourselves if it is possible to 
obtain a version of Chevalley's theorem
in the noncommutative setting.

The first step in answering this question is to determine
what is meant by the coinvariants 
in noncommutative variables. In the
commutative case there are two characterizations of this
$\s_n$-module. First, they can be defined as the solution space
of the system of equations obtained by looking at symmetric
functions without constant term as differential operators (e.g. 
$p_2(\haus) f = \nabla f = 0$). The solution space is called
the harmonics of the symmetric group. The
coinvariants can also be defined as the quotient of
$\bbQ[X] / \SymIdeal$, where $\SymIdeal$ is the ideal
generated by all symmetric functions without constant term.
Indeed, in the commutative case these two definitions
lead to isomorphic spaces.

In the noncommutative setting we have to be more careful.
There are several possibilities for the meaning of a
noncommutative derivative. First, we study
study the harmonics of the symmetric
group with regard to the {\sl Hausdorff} derivative \cite{R}, the differential
operator that acts on letters by 
$\haus_a b = \delta_{a,b}$,
and that satisfies Leibniz rule 
$\haus_a (pq) = (\haus_a p) \cdot q + p \cdot (\haus_a q)$.
We define the harmonics of the symmetric group in the 
noncommutative setting as the space of noncommutative polynomial
solutions
of the system of equations obtained by looking at symmetric
functions without constant term as differential operators
with regard to the
Hausdorff derivative.  We denote this space by $\MH_n$.

In section \ref{sec:lie}, we give an elegant characterization of
this space in terms of the free Lie algebra. 
We show that this space satisfies a mixed commutative/noncommutative
version of Chevalley's theorem.  More precisely,
$\polys \simeq \MH_n \otimes \Sym_n$ where $\Sym_n$
is the space of symmetric polynomials in $n$ variables.

In section \ref{sec:coinvariants}, we look at the 
coinvariants of the symmetric group
in noncommutative variables, defined as
the left quotient
\[
\twistCI
\]
where $\NCSIdeal$ is the left ideal generated by the symmetric
functions in $\NCSI_n$ without constant term. 
We obtain the Hilbert series of this space in terms of the number of Wolf's
irreducible generators \cite{W} which we present 
in a precise combinatorial manner in section \ref{sec:Wolfsummary}.

In addition, we show that Chevalley's theorem
holds in the noncommutative setting. More explicitely, we show that
\[
\polys \simeq \NCSI_n \otimes  \twistCI .
\]
This is done by observing that the coinvariants
of the symmetric group that we just described are isomorphic
to the space of harmonic polynomials with respect to
the twisted derivative defined by
$\twist_a (b v) = \delta_{a,b} v$.

\end{section}

\begin{section}{Combinatorics of set partitions} \label{sec:combnotation}
A set partition $A$ of $m$ is a collection of nonempty subsets $A_1, A_2, \ldots, A_k \subseteq
[m] = \{1,2, \ldots, m \}    $  such  that $A_i \cap A_j = \{ \}$ for $i \neq j$ and
$A_1 \cup A_2 \cup \cdots \cup A_k = [ m ]$.  We will indicate that $A$ is a
set partition of $m$ by the notation $A \vdash [m]$.  The subsets $A_i$ are called the parts
of the set partition and the number of nonempty parts is referred to as the length and will be
denoted by $\ell(A)$.

There is a natural mapping from set partitions to integer partitions 
given by $\la(A)=( |A_1|, |A_2|, \ldots, |A_k| )$, where we assume that 
the blocks of the set partition have been listed in weakly decreasing 
order of size.  If $\la$ is a partition of $n$ (integer partition),
we shall use $\ell(\la)$ to refer to the length (the number of parts) of the partition and
$|\la|$ will be the size of the partition (the sum of the parts), while $n_i(\la)$ shall refer
to the number of parts of the partition of size $i$.  As a convention, 
lowercase Greek letters $\la, \mu$ and $\nu$ will be used to represent integer
partitions while uppercase letters $A, B$ and $C$ will be used for set
partitions.

When writing examples of set partitions we will use the notation that the
sets of numbers are separated by the symbols $\sep$ and the entire set
partition is enclosed by $\{$ and $\}$.  For example, $\{ \{1,3,5\}, \{2\}, \{4\} \}$
will be represented in our notation by $\{ 135\vbar 2\vbar4\}$.  Although there is no
order on the parts of a set partition, we will impose an implied order such
that the parts are arranged by increasing value of the smallest element in the
subset.  This implied order will allow us to reference the $i^{th}$ block of the set
partition without ambiguity.

The number of set partitions is well known and given by the 
Bell numbers.  These can be defined
by the recurrence $B_0 = 1$ and $B_n = \sum_{i=0}^{n-1} \binom{n-1}{ i} B_i$. 
The next seven Bell numbers are $1, 2, 5, 15, 52, 203, 877.$

For a set $S = \{ s_1, s_2, \ldots, s_k \}$ 
of integers $s_i$ and an integer
$n$ we will use the notation $S+n$ to represent the set $\{ s_1+n, s_2+n, \ldots, s_k+n\}$.
For $A \vdash [m]$ and $B\vdash [r]$ set partitions with parts $A_i$, $1\leq i \leq \ell(A)$ and
$B_i$, $1 \leq i \leq \ell(B)$ respectively, we will set $A \attach B = \{ A_1, A_2, \ldots,
A_{\ell(A)}, B_1 + m, B_2+m, \ldots, B_{\ell(B)}+m \}$, therefore $A \attach B \vdash
[m+r]$ and that this operation is noncommutatative in the sense that, in general,
$A \attach B \neq B \attach A$ .

There is a natural lattice structure on the set partitions.  We will define
for $A, B \vdash [n]$ that $A \leq B$ if for each $A_i \in A$ there is a
$B_j \in B$ such that $A_i \subseteq B_j$ (otherwise stated, that $A$
is finer than $B$).  The set of set partitions of $[n]$ with this order forms
a poset with rank function given by $n - k$ where $k$ the length of the set partition.
This poset has minimal element $\{1\vbar2\vbar \cdots \vbar n\}$ and maximal element
$\{12\cdots n\}$.  The largest element smaller
than both $A$ and $B$ will
be denoted $A \wedge B = \{ A_i \cap B_j : 1\leq i \leq \ell(A), 1 \leq j \leq \ell(B) \}$
while the smallest element larger than $A$ and $B$ is denoted $A \vee B$.

\begin{exem}
Let $A = \{ 138\vbar24\vbar5\vbar67 \}$ and $B = \{1\vbar238\vbar4567 \}$.  $A$ and $B$
are not comparable in the inclusion order on set partitions.  We
calculate that $A \wedge B = \{ 1\vbar2\vbar38\vbar4\vbar5\vbar67 \}$ and
$A \vee B = \{ 12345678 \}$.
\end{exem}

When a collection of disjoint sets of positive integers is not a set partition because
the union of the parts is not $[n]$ for some $n$, we may lower the values in
the sets so that they keep their relative values so that the resulting collection
is a set partition.  This operation is
referred to as the `standardization' of a set of disjoint sets $A$ and the
resulting set partition will be denoted $st(A)$.

Now for $A \vdash [m]$ and $S \subseteq \{ 1, 2, \ldots, \ell(A)\}$ with 
$S= \{ s_1, s_2, \ldots, s_k\}$, we define
$A_S = st(\{ A_{s_1}, A_{s_2}, \ldots, A_{s_k} \})$
which will be a set partition of $|A_{s_1}| + |A_{s_2}| +
\ldots + |A_{s_k}|$.  By convention $A_{\{\}}$ is the empty set partition.

\begin{exem}
 If $A = \{ 1368\vbar2\vbar4\vbar579 \}$, then\allowdisplaybreaks{
\begin{align*} 
 &A_{\{1\}} = \{ 1234 \}
 &&A_{\{2\}} = \{ 1 \}
 &&A_{\{3\}} = \{ 1 \}\\
 &A_{\{4\}} = \{ 123 \}
 &&A_{\{1,2\}} = \{ 1345\vbar2 \}
 &&A_{\{1,3\}} = \{ 1245\vbar3 \}\\
 &A_{\{1,4\}} = \{ 1246\vbar357 \}
 &&A_{\{2,3\}} = \{ 1\vbar2 \}
 &&A_{\{2,4\}} = \{ 1\vbar234 \}\\
 &A_{\{3,4\}} = \{ 1\vbar234 \}
 &&A_{\{1,2,3\}} = \{ 1356\vbar2\vbar4 \}
 &&A_{\{1,2,4\}} = \{ 1357\vbar2\vbar468 \}\\
 &A_{\{1,3,4\}} = \{ 1257\vbar3\vbar468 \}
 &&A_{\{2,3,4\}} = \{ 1\vbar2\vbar345 \}
 &&A_{\{1,2,3,4\}} = \{ 1368\vbar2\vbar4\vbar579 \}
 \end{align*}}
\end{exem}
 \end{section}

\begin{section}{The Hopf algebra of noncommutative symmetric functions} \label{sec:Hopfalg}

Consider for a fixed $n >0$ the space $\polys$
consisting of the linear span of monomials in noncommuting variables 
$X_n = \{x_1, x_2, \ldots, x_n\}$.  There is a natural $\s_n$ action on the basis elements
defined by
\begin{equation}
\sigma( x_{i_1} x_{i_2} \cdots x_{i_k} ) = x_{\sigma(i_1)} x_{\sigma(i_2)} \cdots x_{\sigma(i_k)}.
\end{equation}
We can therefore consider $\polys$ as both an $\s_n$ module and an algebra
where the product of two monomials is given by the concatenation of the words.

Let $x_{i_1} x_{i_2} \cdots x_{i_m}$ be a monomial in the space $\polys$.  We
will say that the type of this monomial is the set partition $A \vdash [m]$ with the property
that $i_{a} = i_{b}$ if and only if $a$ and $b$ are in the same block of the set partition.  This
set partition will be denoted as $\nabla(i_1, i_2, \ldots, i_m) = A$.  Notice that the length
of $\nabla(i_1, i_2, \ldots, i_m)$ will be equal to the number of different values which
appear in $(i_1, i_2, \ldots, i_m)$.  

The vector space $\NCSI_n$ will be defined as the linear span of the
elements 
$$\M_{A}[X_n] = \sum_{\nabla(i_1, i_2, \ldots, i_m) 
= A} x_{i_1} x_{i_2} \cdots x_{i_m}$$ 
for $A \vdash [m]$,
where the sum is over all sequences with $1 \leq i_j \leq n$.
For the empty set partition, we define by convention $\M_{\{\}}[X_n] = 1$.
If $\ell(A)>n$ we must have that $\M_A[X_n] = 0$.
Since for any permutation $\sigma \in S_n$,
$\nabla(i_1, i_2, \ldots, i_m) =
\nabla(\sigma(i_1), \sigma(i_2), \ldots, \sigma(i_m))$, we also know
$\sigma \M_A[X_n] = \M_A[X_n]$.

\begin{exem} We list below the monomial NCSFs corresponding to set partitions of size $3$
in a polynomial algebra with 4 variables.
$$\begin{array}{ll}
 \M_{\{123\}}[X_4]    &= x_1 x_1 x_1 + x_2 x_2 x_2 + x_3 x_3 x_3 + x_4 x_4 x_4.\\
 \M_{\{12\vbar3\}}[X_4]   &=  x_1 x_1 x_2 + x_1 x_1 x_3 + x_1 x_1 x_4 + x_2 x_2 x_1 
+ x_2 x_2 x_3 + x_2 x_2 x_4 +\\ 
                     &\phantom{====} x_3 x_3 x_1 + x_3 x_3 x_2 + x_3 x_3 x_4 
+ x_4 x_4 x_1 + x_4^2 x_2 + x_4^2 x_3.\\
 \M_{\{13\vbar2\}}[X_4]   & =  x_1 x_2 x_1 + x_1 x_3 x_1 + x_1 x_4 x_1 + x_2 x_1 x_2. 
                     + x_2 x_3 x_2 + x_2 x_4 x_2 +\\ 
		     &\phantom{====}x_3 x_1 x_3 + x_3 x_2 x_3 + x_3 x_4 x_3+ x_4 x_1 x_4 + x_4 x_2 x_4 + x_4 x_3 x_4.\\
 \M_{\{23\vbar1\}}[X_4]   & =   x_2 x_1 x_1 +  x_3 x_1 x_1 + x_4 x_1  x_1 +  x_1  x_2 x_2. 
+ x_3 x_2 x_2 + x_4 x_2  x_2 +\\  
                     &\phantom{====}x_1 x_3 x_3 +  x_2 x_3 x_3 +  x_4 x_3 x_3+  x_1 x_4  x_4 +  x_2 x_4 x_4 +  x_3 x_4 x_4.\\
 \M_{\{1\vbar2\vbar3\}}[X_4]  &= \sum_{\sigma \in S_4} x_{\sigma(1)} x_{\sigma(2)} x_{\sigma(3)}.
\end{array}
$$
\end{exem}

Now let $\NCSI_n$ be the space of polynomials of $\polys$ which
are invariant under the action of $\s_n$.  For any element $f \in \NCSI_n$, if $\nabla( i_1, i_2,
\ldots, i_k ) = \nabla( j_1, j_2, \ldots, j_k )$ then the coefficient of $x_{i_1} x_{i_2} \cdots x_{i_m}$
in $f$
is equal to the coefficient of $x_{j_1} x_{j_2} \cdots x_{j_k}$ in $f$.  We therefore conclude
that $\{ \M_A[X_n] \}_{\ell(A) \leq n}$ is a basis for $\NCSI_n$.  In addition $\NCSI_n$ has
a ring structure where the product in this ring is defined as the natural extension of the
ring structure on $\polys$.



Our motivation for the following definitions is to extend this algebra
to a Hopf algebra.
Define the vector space $\NCSI^n = {\mathcal L}\{ \M_A \}_{A \vdash [n]}$ where
here we have used $\M_A$ as a symbol representing a basis element for $A$ a set partition.
$\NCSI = \bigoplus_{n\geq 0} \NCSI^n$ is now
the space of noncommutative symmetric functions (as opposed to the space of
noncommutative symmetric polynomials).
The degree of a basis element $\M_A$ is given by $|A|$.
This graded vector space is endowed with a product 
map $\mu : \NCSI^n \otimes \NCSI^m \longrightarrow \NCSI^{m+n}$ which is defined on the basis
elements $\M_A \otimes \M_B$ by
\begin{equation}  
\mu( \M_A \otimes \M_B ) := \sum_{C \vdash [m+n]} \M_{C}
\end{equation}
where the sum is over all set partitions $C$ of $m+n$ such that
$C \wedge (\{1\ldots n\} \attach \{1\ldots m\}) = (A \attach B)$.

This definition is chosen to agree with the product map defined on $\M_A[X_n]$
since we have the following proposition.  

\begin{prop}  \label{prop:varprod} 
Let $A \vdash [r]$ and $B \vdash [m]$, we have
\begin{equation}
\M_A[X_n] \M_B[X_n] = \sum_{C} \M_{C}[X_n]
\end{equation}
where the sum is over all set partitions $C$ of $r+m$ such that
$C \wedge (\{1\ldots r\} \attach \{1\ldots m\}) = (A \attach B)$
with $\ell(C) \leq n$.
\end{prop}

\begin{proof}
The coefficient
of any monomial $x_{i_1} x_{i_2} \cdots x_{i_{m+r}}$ in the expression 
$\M_A[X_n] \M_B[X_n]$ as a product in $\polys$
will have the value
either $1$ or $0$.  As we are working in $\polys$ we must 
have that $\ell(\nabla(i_1, \ldots, i_{r+m})) \leq n$.
We note that the coefficient will be $1$ if and only
if $\nabla(i_1, \ldots, i_r) = A$ and $\nabla( i_{r+1}, \ldots, i_{r+m}) = B$.
This will hold if and only if $\nabla(i_1, \ldots, i_{r+m}) \cap \{ 1, \ldots, r \} = A$
and $\nabla(i_1, \ldots, i_{r+m}) \cap \{ r+1, \ldots, r+m \} = B+r$.  This
is exactly equivalent to the condition that 
$\nabla(i_1, \ldots i_{r+m}) \wedge (\{1\ldots r\} \attach \{1\ldots m\}) = A\attach B$.
\end{proof}

We can conclude that for any $n$, 
the map $\phi_{n}: \NCSI \rightarrow \NCSI_n$ 
is a surjective algebra homomorphism
where $\phi_n$  is
defined as the linear function whose action on the basis is given by
$\phi_n( \M_A ) = \M_A[X_n]$ for $\ell(A) \leq n$ and $\phi_n( \M_A ) = 0$
otherwise.   We state this precisely in the following corollary.

\begin{corollary} The map $\phi_n$ is an algebra morphism.  That is,
\begin{equation}
\phi_n( \mu( \M_A \otimes \M_B ) ) = \phi_n(\M_A) \phi_n(\M_B)
\end{equation}
\end{corollary}

Even though it is defined as an abstract algebra, $\NCSI$ can be realized
as the formal series of bounded degree in an infinite number of variables
which are invariant under all permutations of the indices.  The map
$\phi_n$ is the specialization of this algebra so that the variables $x_{n+1} =
x_{n+2} = x_{n+3} = \cdots = 0$.  In fact we have,
\begin{equation}\label{nvarprince}
\phi_n(F) = 0 \hbox{ for all }n\geq 1\hbox{ if and only if }F = 0.
\end{equation}

The algebra $\NCSI$ was originally considered  by Wolf \cite{W} in extending
 the fundamental theorem of symmetric functions to this algebra and later by
 Bergman and Cohn \cite{BeCo}.  
More recently Rosas and Sagan \cite{RS} considered this space to define natural 
bases which generalize  the bases of the symmetric functions. Our point of
 departure is to consider $\NCSI$  as a Hopf algebra so that we may examine
 it from another perspective.

To this end we define a coproduct map $\Delta : \NCSI^n \longrightarrow \bigoplus_{k=0}^n \NCSI^k \otimes \NCSI^{n-k}$
as
\begin{equation}\label{coproduct}
\Delta( \M_A ) = \sum_{S \subseteq [\ell(A)]} \M_{A_S} \otimes \M_{A_{S^c}}
\end{equation}
where $S^c = [\ell(A)] \backslash S$.

Assume that the $X_n$ and $Y_n$ are two sets of variables which each set is
noncommutative but we have the relations $x_i y_j = y_j x_i$.  Let
$\phi_n^X(\M_A)=\M_A[X_n]$ and $\phi_n^Y( \M_A ) =\M_A[Y_n]$, as before.

\begin{prop} \label{prop:varcoprod} For $F \in \NCSI$, we have
\begin{equation*}
\psi \circ (\phi_n^X \otimes \phi_n^Y) \circ \Delta( F ) = F[X_n,Y_n]
\end{equation*}
where $F[X_n,Y_n]$ represents the noncommutative symmetric polynomial
in $2n$ variables with the additional relations mentioned above and
$\psi( f[X_n,Y_n] \otimes g[X_n, Y_n]) =  f[X_n,Y_n] g[X_n, Y_n]$.  
\end{prop}

\begin{proof}
It suffices to prove this relation for the $\M_A$ basis.  We know then
that
\begin{equation}\label{eq:LHSphi}
\psi \circ (\phi_n^X \otimes \phi_n^Y)\circ  \Delta(\M_A) =
\sum_{S \subseteq [\ell(A)]} \M_{A_S}[X_n]  \M_{A_{S^c}}[Y_n].
\end{equation}

Now in addition 
\begin{equation*}
\M_A[X_n,Y_n] = \sum_{\nabla(i_1, i_2, \ldots, i_r) = A} x_{i_1} x_{i_2} \cdots x_{i_r}
\end{equation*}
where the sum is over the sequences with $1 \leq i_k \leq 2n$ and
we are identifying $x_{i+n} = y_i$ for $1 \leq i \leq n$.  Now for each part
of $A$, $A_i = \{ k_1, k_2, \ldots, k_{|A_i|} \}$, has $i_{k_1} = i_{k_2} = \cdots =
i_{k_{|A_i|}}$.  For a fixed $S \subseteq [\ell(A)]$, consider only the terms
with the property that 
if $i \in S$ and $i_k \in A_i$ then $1 \leq i_k \leq n$ and if $i \notin S$ then
all $i_k \in A_i$ will have $n+1 \leq i_k \leq 2n$ (that is $x_{i_k} = y_{i_k-n}$).
If we restrict the sum to these sequences, then we have
\begin{equation*}
\sum_{\substack{\nabla(i_1, i_2, \ldots, i_r) = A\\
i_k<n+1 \iff i_k \in A_i, i \in S}} x_{i_1} x_{i_2} \cdots x_{i_r} = \M_{A_S}[X_n]
\M_{A_{S^c}}[Y_n].
\end{equation*}
This implies 
\begin{align*}
\M_A[X_n,Y_n] &= \sum_{S \subseteq [\ell(A)]} 
\sum_{\substack{\nabla(i_1, i_2, \ldots, i_r) = A\\
i_k<n+1 \iff i_k \in A_i, i \in S}} x_{i_1} x_{i_2} \cdots x_{i_r}\\
&= \sum_{S \subseteq [\ell(A)]} \M_{A_S}[X_n] \M_{A_{S^c}}[Y_n]
\end{align*}
and this is equal to \eqref{eq:LHSphi}.
\end{proof}
%

In order to have a Hopf algebra we need, in addition, that
the coproduct map is an algebra morphism in the following sense.

\begin{prop} \label{deltacircmu} Let $\tau(F \otimes G) = G \otimes F$ for $F,G \in \NCSI$, then
\begin{equation} \label{eq:Hopf}
\Delta \circ \mu = (\mu \otimes \mu) \circ (id \otimes \tau \otimes id) 
\circ (\Delta \otimes \Delta).
\end{equation}
\end{prop}

\begin{proof}
We will use the previous results and \eqref{nvarprince} to derive this identity.
First we note that  for $F, G \in \NCSI$, and for any $n$ we have by Proposition
\ref{prop:varprod} and Proposition \ref{prop:varcoprod},
$$\psi \circ (\phi_n^{X} \otimes \phi_n^{Y}) \circ \Delta \circ \mu (F \otimes G)
= F[X_n, Y_n] G[X_n, Y_n].$$
The fact that the $X_n$ and $Y_n$ variables commute implies that $\psi \circ \tau \circ (\phi_n^{X} \otimes \phi_n^{Y}) = \psi \circ (\phi_n^{X} \otimes \phi_n^{Y})$.

Therefore, since 
\begin{align*}
\psi \circ (\phi_n^{X} &\otimes \phi_n^{Y})\circ (\mu \otimes \mu) \circ (id \otimes \tau \otimes id)
\circ (\Delta \otimes \Delta)(\M_A \otimes \M_B)\\
&=  \sum_{S \subseteq [\ell (A)]} \sum_{T \subseteq [\ell (B)]} 
 \M_{A_S}[X_n]  \M_{B_T}[X_n]    \M_{A_{S^c}}[Y_n]  \M_{B_{T^c}}[Y_n]   \\
 &=  \sum_{S \subseteq [\ell (A)]} 
 \M_{A_S}[X_n]     \M_{A_{S^c}}[Y_n]    
 \sum_{T \subseteq [\ell (B)]}    \M_{B_T}[X_n]  \M_{B_{T^c}}[Y_n]   \\
&=\M_A [X_n, Y_n] \M_B [X_n, Y_n].
\end{align*}



Now, since it suffices to prove the relation for the $\M_A$ basis and by Proposition \ref{prop:varcoprod} 
$$\psi \circ (\phi_n^{X} \otimes \phi_n^{Y}) \circ \Delta \circ \mu 
=\psi \circ (\phi_n^{X} \otimes \phi_n^{Y})\circ (\mu \otimes \mu) \circ (id \otimes \tau \otimes id)
\circ (\Delta \otimes \Delta)
$$
 holds for any $n$, we must have that 
equation (\ref{eq:Hopf}) holds on $\NCSI$.
\end{proof}
\end{section}

\begin{section}{Hopf Algebras}
One of the main reasons for looking at this space as a Hopf
algebra is that we are able to put it in context with other
well known Hopf algebras.  To this end, we introduce the
space of symmetric functions $\Sym$ and another algebra referred to
as the noncommutative symmetric functions $\NSym$.

For each graded bialgebra $\H$,
we will have an implicit unit map $u^{\H}$ sending the $1$ in the field
to the degree $0$ basis element (also denoted by $1$) and the
counit $\varepsilon^{\H}$ which sends all terms of degree greater
than $0$ to $0$ and the $1$ of the algebra to the $1$ of our
base field (in the following algebras will always be $\bbQ$).

In each of the following bialgebras, the product $\mu^{\H}$ and
coproduct $\Delta^{\H}$ respect the grading in the sense that
$\mu^{\H} : \H^n \otimes \H^m \rightarrow \H^{n+m}$ and
$\Delta^{\H} : \H^n \rightarrow \bigoplus_{k=0}^n \H^{k} \otimes
\H^{n-k}$ where we will indicate the vector spaces $\H^n$
as the homogeneous components of the algebra of degree $n$.
It is a well known result that every graded bialgebra where
the degree $0$ component has dimension $1$ is a Hopf algebra
(i.e., is a connected Hopf algebra)
\cite{Sw}.

We will not give the antipode explicitly on the algebras but
it is defined uniquely for any graded, connected Hopf algebra
by the defining relation:
\begin{equation*}
\mu^{\H}\circ(id \otimes S^{\H}) \circ \Delta^{\H} = u^{\H}\circ\varepsilon^{\H}.
\end{equation*}
To compute the action of the antipode on an element $f$ of homogeneous
degree greater than $0$, we write $\Delta^{\H}(f) = 1 \otimes f + \sum_{i} g_{i}' \otimes
g_{i}''$.  It then follows that
\begin{equation*}
0 = \mu^{\H}\circ(id \otimes S^{\H}) \circ \Delta^{\H}(f) 
= S^{\H} (f) + \sum_{i} \mu^{\H} (g_{i}' \otimes
S^{\H}(g_{i}''))
\end{equation*}
which can be used to solve for $S^{\H} (f)$ while the $S^{\H}(g_{i}'')$ will be
of smaller degree and can be computed recursively.

From this discussion and Proposition \ref{deltacircmu} we can conclude
\begin{theorem}
$\NCSI$ is a Hopf algebra.
\end{theorem}

\begin{subsection}{Symmetric functions}
We will define the algebra of symmetric functions $\Sym$ as the 
free commutative algebra generated by elements $h_k$ for $k \geq 1$.
The product on this algebra is the standard commutative product with
a grading defined by $deg(h_k) =k$ and following convention we will
denote $\Sym^n = {\mathcal L}\{ h_\la \}_{\la \vdash n}$ with $h_\la :=
h_{\la_1} h_{\la_2} \cdots h_{\la_{\ell(\la)}}$
and set $\Sym = \bigoplus_{n \geq 0} \Sym^n
= \bbQ[ h_1, h_2, h_3, \ldots]$.  

We may define the graded dual algebra $\Sym^*$ by defining
the graded component of degree $n$ is the vector space defined
$(\Sym^\ast)^n = {\mathcal L}\{ m_\la \}_{\la \vdash n}$  and 
space is $\Sym^\ast = \bigoplus_{n \geq 0} (\Sym^\ast)^n$ where the
basis $m_\la$ is dual to $h_\la$ in the dual pairing.  It arises that
$\Sym^\ast \simeq \Sym$. In fact,
\begin{equation}
h_n = \sum_{\la \vdash n} m_\la.
\end{equation}
The product and coproduct on the $m_\la$ basis can be determined
from the product and coproduct on the $h$-basis.  It develops that,
\begin{equation}\label{eq:mprod}
m_\la m_\mu = \sum_{\nu} r^{\nu}_{\la\mu} m_\nu
\end{equation}
where the coefficients $r^\nu_{\la\mu}$ are the number of pairs of vectors
$(\alpha, \beta)$ such that $\alpha \sim \la$, $\beta \sim \mu$ such that
$\alpha_i + \beta_i = \nu_i$ for all $i$ and $\alpha \sim \la$ means that the sequence of
values of $\alpha$ rearranges to the partition $\la$.  

The coproduct is given by the formulas
\begin{equation}
\Delta^\Sym( h_n ) = \sum_{k=0}^n h_k \otimes h_{n-k}
\end{equation}
and
\begin{equation}\label{eq:mcoprod}
\Delta^{\Sym^*}( m_\la ) = 
\sum_{\mu\uplus\nu=\la} m_\mu \otimes m_\nu
\end{equation}
where we have denoted in the sum by $\mu\uplus\nu=\la$
that $\mu$ and $\nu$ are partitions satisfying
$n_i(\nu) + n_i(\mu) = n_i(\la)$ for all $i\geq 1$.
\end{subsection}

\begin{subsection}{The other Hopf algebra of noncommutative symmetric functions}
\label{sec:HopfNSym}
There exists a noncommutative algebra which can be seen as an
analogue to $\Sym$ (while 
$\NCSI = \bigoplus_{n \geq 0} {\mathcal L}\{ \M_A \}_{A \vdash [n]}$ 
is more clearly an analogue of 
$\Sym^\ast= \bigoplus_{n \geq 0} {\mathcal L}\{ m_\la \}_{\la \vdash n}$), see 
for example \cite{thibon}.  
$\NSym$ is defined as the noncommutative polynomial ring generated freely by
elements $\h_k$ for $k \geq 1$ where $deg(\h_k) = k$.  For a composition
$\alpha$ of $n$ (denoted $\alpha \models n$)
we set $\h_\alpha := \h_{\alpha_1} \h_{\alpha_2} \cdots \h_{\alpha_{\ell(\alpha)}}$
and $\NSym^n = {\mathcal L}\{ \h_{\alpha} \}_{\alpha \models n}$.  
The graded algebra is then defined as
$$\NSym = \bigoplus_{n \geq 0} \NSym^n = \bbQ\langle \h_1, \h_2, \h_3, \ldots \rangle.$$
The compositions of $n$ are in bijection with the subsets of $\{ 1, 2, \ldots, n-1\}$
by the correspondence $D(\alpha) = \{ \alpha_1, \alpha_1+\alpha_2, \ldots, \alpha_1+
\alpha_2+\cdots+\alpha_{\ell(\alpha)-1}\}$ (the descent set of
the composition $\alpha$) and hence $\dim( \NSym^n ) = 2^{n-1}$.

The product on $\NSym$ is defined so that 
$\mu^\NSym : \NSym^n \otimes \NSym^m \rightarrow \NSym^{n+m}$
as the free noncommutative product $\mu^\NSym( \h_\alpha \otimes \h_\beta ) =
\h_\alpha \h_\beta$.

The coproduct is given by the following formula and the fact that $\Delta^\NSym$
is an algebra homomorphism.
\begin{equation}
\Delta^\NSym( \h_n ) = \sum_{k=0}^n \h_k \otimes \h_{n-k}
\end{equation}

There is a significant difference between the dimensions of the two
Hopf algebras of noncommutative symmetric functions. The dimension
of $\NCSI^n$   is the number of set partitions of $[n]$,  and for $n>2$
this is larger than the dimension of $\NSym^n$ which is the number
of compositions of $n$ or
$\dim(\NSym^n) = 2^{n-1}$.
\end{subsection}

\begin{subsection}{Relations between $\Sym$, $\NSym$, and $\NCSI$}
The symbol $\chi$ will represent the `forgetful' map which sends
elements of a noncommutative algebra to the commutative counterpart
(the map which `forgets' that expressions are noncommutative).
In our case we will begin by considering two such maps.  The
first of which is $\chi : \NSym \rightarrow \Sym$ given by the linear
homomorphism $\chi( \h_\alpha) = h_{\alpha_1} h_{\alpha_2} \cdots h_{\alpha_{\ell(\alpha)}}$.

\begin{prop}  The linear map
$\chi : \NSym \rightarrow \Sym$ where $\chi( \h_\alpha) = h_\alpha$
is a Hopf morphism.
\end{prop}
\begin{proof}
This is easy to check on the $\h_\alpha$ basis since 
$\chi(\h_{\alpha} \h_{\beta})
= h_\alpha h_\beta$ and $$(\chi\otimes\chi)\circ\Delta^\NSym( \h_n ) = \sum_{k=0}^n h_k \otimes
h_{n-k} = \Delta^\Sym \circ \chi(\h_n).$$  Since both $\chi$ and
$\Delta$ are algebra homomorphisms,
this relation will hold as well on basis elements $\h_\alpha$.
\end{proof}


In addition we will use the same symbol $\chi$ to represent the
map $\chi: \NCSI \rightarrow \Sym \cong \Sym^\ast$ given by the linear homomorphism
$\chi( \M_A ) = \la(A)^! m_{\la(A)}$ where we denote
$\la^! = (\prod_{i\geq1} n_i(\la)!)$. By contrast, we will use $\la! = \la_1! \la_2! \cdots
\la_{\ell(\la)}!$ (these conventions use
the notation introduced in section \ref{sec:combnotation} 
and are consistent with the notation of \cite{RS}).  This map is inspired
from the expression $\M_A[X_n]$ since if the variables were 
allowed to commute then the expression is equal to $\la(A)^! m_{\la(A)}[X_{n}]$
where for a partition $\la$, $m_{\la}[X_{n}] = 
\sum_{\alpha \sim \la} x^{\alpha}$ is the monomial symmetric polynomial.

\begin{prop}\label{prop:chiisHopf}  The linear map 
$\chi : \NCSI \rightarrow \Sym$ where 
$\chi( \M_A) = \la(A)^!\,m_{\la(A)}$ is a Hopf morphism.
\end{prop}

\begin{proof}
As we remarked above, $\chi: \NCSI \rightarrow \Sym$ is the restriction
of the map
\begin{equation*}
\chi: \bbQ\langle x_{1}, x_{2}, x_{3}, \ldots \rangle \longrightarrow
\bbQ[ x_{1}, x_{2}, x_{3}, \ldots ]
\end{equation*}
which is the map that forgets the variables are noncommutative.  Clearly
this map is an algebra morphism.  It follows that the restriction of
this map to $\NCSI$ and $\Sym$ will also be an algebra morphism.

In addition we need to show that $(\chi \otimes \chi) \circ \Delta^\NCSI =
\Delta^\Sym \circ \chi$.
We remark that for a given $\mu$ the number of 
subsets $S$ such that $\la(A_S)=\mu$ and
$\la(A_{S^c}) = \nu$ with $n_i(\mu) + n_i(\nu) = n_i(\la(A))$ is
equal to $\frac{\la(A)^!}{\mu^! \nu^!}$.  Therefore
\begin{align*}
(\chi \otimes \chi) \circ \Delta^\NCSI( \M_A) &=
\sum_{S \subseteq [\ell(A)]} \la(A_S)^!\,\la(A_{S^c})^!\,m_{\la(A_S)} \otimes m_{\la(A_{S^c})}\\
&= \sum_{\mu\uplus\nu=\la(A)} \mu^! \nu^!
\left( \frac{\la(A)^!}{\mu^! \nu^!}\right)
m_{\mu} \otimes m_{\nu}\\
&= \sum_{\mu\uplus\nu=\la(A)} \la(A)^!
m_{\mu} \otimes m_{\nu}\\
&= \Delta^\Sym(\la(A)^! m_{\la(A)}) = \Delta^\Sym( \chi(\M_A)).
\end{align*}
Therefore $\chi$ is also a morphism with respect to the coproduct and
hence is a Hopf morphism.
\end{proof}

There is a natural pullback of $\chi:\NCSI \rightarrow \Sym$ which
was considered by Rosas and Sagan in \cite{RS}. 
They called the linear homomorphism
${\tilde \chi} : \Sym \rightarrow \NCSI$ defined by ${\tilde \chi}(m_\la) = \frac{\la!}{|\la|!}
\sum_{\la(A) = \la} \M_A$ the {\sl lifting map} and showed it has the following
property.
\begin{prop} \cite[Proposition 4.1]{RS}
$\chi\circ{\tilde \chi}$ is the identity map on $\Sym$.
\end{prop}

For our purposes, the important property of the lifting map
will be from the following proposition.
\begin{prop} \label{prop:liftcomorph}
\begin{equation}
\Delta^\NCSI \circ {\tilde \chi} = ({\tilde \chi} \otimes {\tilde \chi}) \circ
\Delta^\Sym.
\end{equation}
\end{prop}
\begin{proof}
 Equation (\ref{coproduct}) allows us to  deduce this property  by direct computation.
\begin{align*}
\Delta^\NCSI ( {\tilde \chi}( m_\la))&= \frac{\la!}{|\la|!} \sum_{C: \la(C) = \la}
\Delta^\NCSI( \M_C)\\
&= \frac{\la!}{|\la|!} \sum_{C: \la(C) = \la}
\sum_{S \subseteq [\ell(\la)]} \M_{C_S} \otimes \M_{C_{S^c}}\\
&= \frac{\la!}{|\la|!} \sum_{C: \la(C) = \la}
\sum_{\mu\uplus\nu=\la}
\sum_{\substack{S \subseteq [\ell(\la)]\\
\la(C_S)=\mu}} \M_{C_S} \otimes \M_{C_{S^c}}\\
&= \frac{\la!}{|\la|!} \sum_{C: \la(C) = \la}
\sum_{\mu\uplus\nu=\la} \sum_{\substack{A:\la(A) = \mu\\
B:\la(B)=\nu}}
\sum_{\substack{S \subseteq [\ell(\la)]\\
C_S=A\\
C_{S^c}=B}} \M_{A} \otimes \M_{B}
\end{align*}
Now to complete this computation we exchange the sums and notice that for
a fixed set partitions $A$ and $B$, there are exactly $\binom{|C|}{ |A|}$
different set partitions $C$
such that there is an $S \subseteq [\ell(C)]$ with $C_S = A$ and $C_{S^c} = B$.
\begin{align*}
&= \frac{\la!}{|\la|!} 
\sum_{\mu\uplus\nu=\la} \sum_{\substack{A:\la(A) = \mu\\
B:\la(B)=\nu}} \binom{|\la| }{ |\mu|} \M_{A} \otimes \M_{B}\\
&=  
\sum_{\mu\uplus\nu=\la} \sum_{\substack{A:\la(A) = \mu\\
B:\la(B)=\nu}} \frac{\la!}{|\la|!} \frac{|\la|!}{|\mu|!|\nu|!} \M_{A} \otimes \M_{B}\\
&=  
\sum_{\mu\uplus\nu=\la} \sum_{\substack{A:\la(A) = \mu\\
B:\la(B)=\nu}} \frac{\la!}{|\mu|!|\nu|!} \M_{A} \otimes \M_{B}.
\end{align*}
Now since $\mu\uplus\nu=\la$ then we have that $\la!=\mu!\nu!$ and hence the
equation above is equal to
\begin{equation*}
=\sum_{\mu\uplus\nu=\la} {\tilde \chi}(m_\mu) \otimes {\tilde \chi}(m_\nu)
= ({\tilde \chi}\otimes{\tilde \chi})\circ\Delta^\Sym(m_\la)
\end{equation*}
\end{proof}

This last proposition leads us to identifying an important relationship between
the algebra of $\NSym$ of noncommutative symmetric functions and the
algebra of $\NCSI$.

\begin{theorem} Define $\NSaNC : \NSym \rightarrow \NCSI$ by the action 
$$\NSaNC( \h_n ) = {\tilde \chi}( h_{n}) = \sum_{A \vdash [n]} \frac{\la(A)!}{n!} \M_A$$
and extend this as an algebra morphism by defining for the linear basis $\h_\alpha$,
\begin{equation*}
\NSaNC( \h_\alpha ) = {\tilde \chi}( h_{\alpha_1}) 
{\tilde \chi}( h_{\alpha_2}) \cdots 
{\tilde \chi}( h_{\alpha_{\ell(\alpha)}}).
\end{equation*}
$\NSaNC$ is a Hopf morphism and $\NSaNC$ an inclusion map so
that $\NSym$ is a natural subalgebra of
$\NCSI$.
\end{theorem}

Before proceeding with the proof of the theorem we introduce
an important lemma.
For each $\alpha$ a composition of $n$, we have a canonical corresponding set
partition, 
\begin{equation}\label{eq:compsp}
A(\alpha)=\{1,2,\ldots,\alpha_1\vbar \alpha_1+1,\ldots,\alpha_1+\alpha_2 \vbar
\cdots \vbar \alpha_1+\cdots+\alpha_{\ell(\alpha)-1},\ldots,|\alpha|\}.
\end{equation}

\begin{lemma} \label{lem:coeff}
The coefficient of
$\M_{A(\beta)}$ in $\NSaNC(\h_\alpha)$ is equal to $(\alpha \cup \beta)!/\alpha!$
where $\alpha \cup \beta$ is the composition with descent set equal to
$D(\alpha) \cup D(\beta)$ and  $\alpha! = \alpha_1! \alpha_2! \cdots \alpha_{\ell(\alpha)}!$
for a composition $\alpha$.
\end{lemma}

\begin{proof}
Note that $\NSaNC(\h_\alpha) = \NSaNC( \h_{\tilde \alpha} ) \NSaNC(\h_{\alpha_{\ell(\alpha)}})$
where ${\tilde \alpha} = (\alpha_1, \alpha_2, \ldots, \alpha_{\ell(\alpha)-1})$.
Let ${\tilde B} = A(\beta)\coeff_{\{1,\ldots,|\tilde \alpha|\}}$ which corresponds to a
composition ${\tilde \beta}$ and ${\bar B} = st( A(\beta) 
\coeff_{\{|\tilde \alpha|+1,\ldots,|\alpha|\})})$.
When $\M_{A(\beta)}$ arises as the coefficient of $\NSaNC(\h_\alpha)$ the coefficient will
be the coefficient of $\M_{\tilde B}$ in $\NSaNC( \h_{\tilde \alpha} )$
times the coefficient of $\M_{\bar B}$ in $\NSaNC( \h_{\alpha_{\ell(\alpha)}})$.  By induction
on the number of parts of $\alpha$ we can assume that this coefficient
is $\frac{({\tilde \alpha} \cup {\tilde \beta})!}{{\tilde \alpha}!} 
\frac{\lambda({\bar B})!}{\alpha_{\ell(\alpha)}!} = \frac{({\tilde \alpha} \cup {\tilde \beta})!
\lambda({\bar B})!}{\alpha!} = \frac{(\alpha \cup \beta)!}{\alpha!}.$
\end{proof}

\begin{proof}[Proof of Theorem]
Proposition \ref{prop:liftcomorph} says that
\begin{align*}
\Delta^\NCSI(\NSaNC( \h_n)) &= ({\tilde \chi} \otimes {\tilde \chi}) \circ
\Delta^\Sym(h_n) \\
&= \sum_{k=0}^n {\tilde \chi}(h_k) \otimes {\tilde \chi}(h_{n-k}) \\
&= (\NSaNC \otimes \NSaNC) \circ \Delta^\NSym( \h_n).
\end{align*}
Clearly we have that $\NSaNC(\h_\alpha \h_\beta) = \NSaNC(\h_\alpha)
\NSaNC(\h_\beta)$ so we know that $\NSaNC$ is a Hopf morphism.  In order to
show that
$\NSaNC$ is an inclusion of $\NSym$ into $\NCSI$ we need to show that
the generators of $\NSym$, $\NSaNC(\h_n)$, are algebraically independent. 
This is equivalent to showing
that the elements $\NSaNC(\h_{\alpha_1}) \NSaNC(\h_{\alpha_2})\cdots
\NSaNC(\h_{\alpha_{\ell(\alpha)}})$ are linearly independent.

 In order to show that $\NSaNC(\h_\alpha)$
are linearly independent, it suffices to examine
the minor of coefficients of $\M_{A(\beta)}$ in $\NSaNC(\h_\alpha)$ and show
that the determinant of this minor is nonzero.  The coefficient of
 $\M_{A(\beta)}$ in $\NSaNC(\h_\alpha)$ is $(\alpha \cup \beta)!/\alpha!$ by
 Lemma \ref{lem:coeff}.

The proof follows by showing that the $2^{n-1} \times 2^{n-1}$ determinant
of the matrix $\left[ (\alpha \cup \beta)! \right]_{\alpha,\beta \models n}$
is nonzero.  In Theorem \ref{thrm:magicdet} below we compute that
this matrix has a nonzero determinant (in fact we compute it explicitly) and hence
conclude that $\NSaNC$ is an inclusion.
\end{proof}

By writing the first few matrices and their determinants gives
a clue on how to show that it has a nonzero determinant.
Begin by ordering the compositions in lexicographic order so that
$(11) < (2)$, $(111)<(12)<(21)<(3)$, and $(1111)<(112)<(121)<(13)<(211)<(22)<(31)<(4)$
are the order of the indices of the matrices below.

\begin{equation*}
\begin{array}{cc}
\vline~1&1~\vline\\
\vline~1&2~\vline
\end{array} = 1
\end{equation*}

\begin{equation*}
\begin{array}{cccc}
\vline~1&1&1&1~\vline\\
\vline~1&2&1&2~\vline\\
\vline~1&1&2&2~\vline\\
\vline~1&2&2&6~\vline\\
\end{array} = 3
\end{equation*}
\begin{equation*}
\begin{array}{cccccccc}
\vline~1&1&1&1&1&1&1&1\phantom{1}~\vline\\
\vline~1&2&1&2&1&2&1&2\phantom{1}~\vline\\
\vline~1&1&2&2&1&1&2&2\phantom{1}~\vline\\
\vline~1&2&2&6&1&2&2&6\phantom{1}~\vline\\
\vline~1&1&1&1&2&2&2&2\phantom{1}~\vline\\
\vline~1&2&1&2&2&4&2&4\phantom{1}~\vline\\
\vline~1&1&2&2&2&2&6&6\phantom{1}~\vline\\
\vline~1&2&2&6&2&4&6&24~\vline\\
\end{array} = 117 = 3^2\cdot13
\end{equation*}
The next two values of this determinant
are $2915757 = 3^5\cdot13^2\cdot71$ and $458552896435013913 = 3^{12} \cdot13^5\cdot71^2\cdot461$.  Although the sequence of determinants
is not familiar, the factors which appear in it are.  The sequence
$1, 1, 3, 13, 71, 461,\ldots$ are found in the OLEIS \cite{Sl} as sequence
A003319, the permutations of $n$ with no global descents.  A global
descent is a value $k$ such that $\pi_i > \pi_j$ for all $i \leq k$ and $j>k$.  The
number of these can be calculated with the recurrence $a_{1}=1$ and for $n>1$,
\begin{equation} \label{eq:anformula}
a_{n} = n! - \sum_{i=1}^{n-1} a_{i} (n-i)!.
\end{equation}  
The permutations with no global descents
arise in the Hopf algebra of permutations
due to Malvenuto-Reutenauer as the primitive elements/generators
of the Hopf algebra \cite{AS}.  This begs an explanation of why these numbers
should arise in this computation.
The expression for the determinant
is summarized in the following theorem.

\begin{theorem} \label{thrm:magicdet}
\begin{equation} \label{eq:magicdet}
det\left| (\alpha \cup \beta)! \right|_{\alpha,\beta \models n} = \prod_{\alpha \models n}
\prod_{i=1}^{\ell(\alpha)} a_{\alpha_i}
\end{equation}
where $a_n$ is the number of permutations of $n$ with no global descents.
\end{theorem}
%


One could prove this theorem by induction using the identity for
a block matrix (see \cite{Mo})
$$ det \left[ \begin{array}{cc}
A&B\\C&D\end{array}\right] = det(A) det( D- C A^{-1} B)
$$
and the recursive structure of the matrix 
$\left[ (\alpha \cup \beta)! \right]_{\alpha, \beta \models n}$.
Instead we present a proof suggested to us by A. Lascoux \cite{L} which
makes the formula of Theorem \ref{thrm:magicdet} transparent.

\begin{proof}
It is easy to see that any permutation can be decomposed  uniquely as a concatenation
of permutations with no global descents. For example, $ 4 6 5 3 1 2$ can be decomposed as 
$ 4 6 5 \cdot 3 \cdot 1 2$.

Using this decomposition, we associate to any permutation a composition. 
For instance to $ 4 6 5 \cdot 3 \cdot 1 2$
we associate the composition $(3, 1, 2)$. Therefore, the set of all permutations 
in $\s_n$ can be partitioned into a disjoint union of subclasses indexed
by compositions.

In addition, the cardinality of the class indexed by
$\alpha \models n$ is $a_\alpha:=a_{\alpha_1} a_{\alpha_2} \cdots a_{\alpha_{\ell(\alpha)}}$, 
where $a_j$ is the number of permutations with no global descents. 

We conclude that
$$\sum_{\alpha \models n} a_\alpha = n!~.
$$

An analogous statement holds for any Young subgroup of $\s_n$. To wit,
\begin{equation}\label{key}
\sum_{\beta \leq \alpha  } a_\beta = \alpha!~
\end{equation}
with $\leq$ representing the standard refinement order on compositions.
Now our notation for $\alpha \cup \beta$ means that for
$\eta \leq (\alpha \cup \beta)$ that
$\eta \leq \alpha$ and $\eta \leq \beta$ and so using the notation
$\dop true \dcp = 1$ and $\dop false \dcp = 0$ we have the expression
\begin{equation}
(\alpha \cup \beta)!  = \sum_{\eta \leq (\alpha \cup \beta) } a_\eta
= \sum_{\eta} \dop \eta \leq \alpha\dcp \dop \eta \leq \beta \dcp \, a_\eta.
\end{equation}

Let $\mathbb D$ represent the $2^{n-1} \times 2^{n-1}$ 
diagonal matrix indexed by compositions $\eta$ with 
$a_\eta$ the entries along the diagonal.  Also let ${\mathbb C} =
[ \dop\beta \leq \alpha\dcp ]_{\alpha, \beta \models n}$, a matrix with entry $1$
at $(\alpha, \beta)$ if $\beta \leq \alpha$ and $0$ otherwise.  
Now look at the $(\alpha, \beta)$ entry in the product
${\mathbb C} {\mathbb D} {\mathbb C}^T$.  This will be
$$\sum_{\delta, \theta \models n} \dop \delta \leq \alpha \dcp \dop \delta = \theta \dcp
\, a_\theta \, \dop \theta \leq \beta\dcp
= \sum_{\delta \models n} \dop \delta \leq \alpha \dcp \dop \delta \leq \beta \dcp  \, a_\delta =
(\alpha \cup \beta)!~.$$
We conclude that ${\mathbb C} {\mathbb D} {\mathbb C}^T = 
\left[ (\alpha \cup \beta)! \right]_{\alpha,\beta \models n}$ and hence
$\det \left[ (\alpha \cup \beta)! \right]_{\alpha,\beta \models n}  = \det  {\mathbb D}$
(since $\det {\mathbb C} = 1$).  This demonstrates \eqref{eq:magicdet} since $\mathbb D$
is a diagonal matrix with determinant equal to $\prod_{\alpha \models n} a_\alpha$.
\end{proof}

\begin{remark}
$\NSym$ is also generated by the analogs of the power and
elementary bases of the symmetric functions and there are
formulas for expressing these into the $\h$-basis.
The map $\NSaNC$ is not unique since we could just as easily lift
these other bases (as we defined $\NSaNC( \h_n ) = {\tilde \chi}( h_n)$)
and by direct computation one can verify that these other inclusions
of $\NSym$ in $\NCSI$ are not the same as $\NSaNC$.
For example, $\e_3 = \h_3 - \h_{12} - \h_{21} + \h_{111}$ and
$$\NSaNC(\e_3) = \frac{1}{6} \M_{\{1\vbar2\vbar3\}}
+\frac{1}{3} \M_{\{13\vbar2\}} - \frac{1}{6} \M_{\{12\vbar3\}} - \frac{1}{6} \M_{\{1\vbar23\}} \neq
{\tilde \chi}( e_3 ) = \frac{1}{6} \M_{\{ 1\vbar2\vbar3 \}}.$$
Therefore the inclusion that we present here is not
unique, but once we fix a set of generators of $\Sym$ there is a natural
embedding of $\NSym$ to $\NCSI$.
\end{remark}


We conclude this section with a summary of these results
stating that the
Hopf algebra morphisms which
relate $\NSym$, $\NCSI$ and $\Sym$ can be drawn in a commutative
diagram.
\begin{theorem}
The following diagram commutes and all maps are Hopf morphisms.
\begin{equation*}
\begin{CD}
\NSym   @>{\NSaNC}>>   \NCSI\\ 
@VV{\chi}V                                    @VV{\chi}V\\ 
\Sym                  @=      \Sym^{\ast}
\end{CD}
\end{equation*}
\end{theorem}
\vskip .2in

\section{Remarks on the algebra structure of noncommutative symmetric functions} \label{sec:algebra}

There are other relationships between $\NCSI$ and $\NSym$ that
are worth considering but are not as structured because they only
hold on the level of algebras and do not respect the coproduct.  
First, we shall examine a graded algebra isomorphism between the 
graded algebras  $\NSym$ and $\NCSI_2$. 

This algebra isomorphism implies that the structure 
constants for $\NCSI_2$ with respect to the monomial basis coincide 
with the structure constants for $\NSym$ in the ribbon 
Schur basis,  which are know to be related to the
representation theory of the Hecke algebra at $q=0$, see
for example \cite{KrTh}.  

 
In general, the structure constants of $\NCSI_n$ with respect to
the monomial basis (as well as those 
of $\NCSI$) are also nonnegative integers. A natural question to 
ask is whether the representation theoretical interpretations of 
$\NCSI_2$ can be extended to $\NCSI$, as well as to its specializations 
$\NCSI_n$, for each value of $n$.

The number of set partitions of $[n]$ with at most  two blocks is $2^{n-1}$. 
Therefore, there is a bijection between  set partitions of $[n]$ with at 
most two parts and compositions of $n$.  Compositions are the indexing
set of the algebra $\NSym$ of subsection \ref{sec:HopfNSym} and so it
is natural to look for a connection through this structure. 
In fact, this observation gives us a 
way of relating $\NCSI_2$ and $\NSym$. That is, $\NSym$ is isomorphic, as 
an algebra, to  $\NCSI_2$.

To any set partition $A= \{A_1\vbar A_2 \}$, we associate the ribbon shape obtained
by reading numbers $1, 2, \cdots, n$ sequentially, and placing the $(i+1)^{st}$
box to the right of the $i^{th}$ box if  $i$ and $i+1$ are in the same block 
of $A$, or placing the  $(i+1)^{st}$ box immediately  below the $i^{th}$ 
box otherwise.
For instance, the ribbon associated to $A = \{1 2 4 5 \vbar 3 6 7 8\}$ is
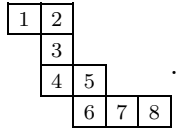
   $$\begin{matrix} 
               \Young{\Box{1} & \Box{2}\cr  &\Box{3}\cr
                & \Box{4} & \Box{5}   \cr  & &\Box{6} & \Box{7} & \Box{8} \cr } 
     \end{matrix} .  $$
Note that we are placing numbers inside the boxes 
for the sake of clarity only. We denote by $c(A)$ the composition of $n$ 
obtained by recording the lengths of the horizontal segments in the 
corresponding ribbon. In our example, $c(A)=(2,1,2,3).$
     
Following \cite{thibon}, we define a second and very important basis for 
$\NSym$, the ribbon Schur functions. Recall that the set of all compositions of 
$n$ is equipped with the reverse refinement order of the
descent sets, denoted  by $\le$. For 
instance, $(2,2,1) \le (4,1)$.     The ribbon Schur functions
$({\bf R}_{\alpha})$ are defined by the following expression:
\begin{align*}
\bf{R}_{\alpha} &= \sum_{\beta \le \alpha} (-1)^{\ell(\alpha) - \ell(\beta)}
\h_{\alpha}
\end{align*}

Let $\iota : \NCSI_2 \mapsto \NSym$ by the linear homomorphism such that
$$\iota( \M_{\{ A\vbar B \} } ) = {\bf R}_{c(\{A\vbar B\})},$$ where $c(\{ A\vbar B \} )$ 
denotes the composition corresponding to $\{ A\vbar B \}$ under the bijection 
just stated.

\begin{prop} The map $\iota : \NCSI_2 \mapsto \NSym$ is an isomorphism of algebras.
\end{prop}

\begin{proof}

The monomials in $\NCSI_2$ multiply according to the following rule
$$\M_{\{A_1\vbar A_2\}} \, \M_{\{B_1\vbar B_2\}} = \M_{\{A_1\cup (B_1+|A|)\vbar A_2 \cup (B_2+|A|) \}} 
+ \M_{\{A_1 \cup (B_2+|A|)\vbar A_2 \cup (B_1+|A|)\}}.$$
For instance, if $A= \{ 1 3 4 6 \vbar 2 5 7 8 \} $ and 
$B = \{ 1 2 \vbar 3 4 5 \}$, then
$B + 8 = \{ 9 \,10 \vbar 11\, 12\, 13 \},$  and
the two terms in the product of $\M_A[x_1, x_2]$ and $\M_B[x_1, x_2]$ are indexed by
$$\{A_1 \cup (B_1+8)\vbar A_2 \cup (B_2+8)\}= \{ 1 3 4 6 9\, 10 \vbar 2 5 7 8 \,11\, 12 \,13 \}$$
and  
$$\{A_1\cup (B_2+8)\vbar A_2 \cup (B_1+8)\} =\{ 1 3 4 6 \,11\, 12\, 13\vbar 2 5 7 8 9\, 10 \}.$$

On the other hand, it is well known that the ribbon Schur functions 
multiply as $\bf{R}_{\alpha}\bf{R}_{\beta}=\bf{R}_{\alpha \rhd \beta} 
+ \bf{R}_{\alpha \cdot \beta},$
where $\alpha \rhd \beta$ is the composition obtained by adding the last 
part of $\alpha$ to the first part of $\beta$, and $\alpha \cdot \beta $ 
is the composition obtained by concatenation. Hence, multiplying  ribbons 
$\bf{R}_{\alpha}$ and $\bf{R}_{\beta}$ is equivalent to placing the 
first box of $\bf{R}_{\beta}$ next to the last box of $\bf{R}_{\alpha}$  
either vertically or horizontally. 

To finish our argument, note that if $\alpha=c(\{A_1\vbar A_2\})$, and 
$\beta=c(\{B_1\vbar B_2\})$. Then, $\alpha \rhd \beta =c(\{A_1\cup (B_1+|A|)\vbar 
A_2\cup (B_2+|A|)\}) $ and 
$\alpha \cdot \beta =c(\{A_1\cup(B_2+|A|)\vbar A_2\cup(B_1+|A|)\}).$ This is best done looking at 
our running example and noticing that joining the last row of $\alpha$ and
the first row of $\beta$ 
corresponds to joining blocks $A_1$ and $B_1$ together, and placing the 
the first row of $\beta$ below the last row of $\alpha$
corresponds to joining blocks $A_1$ and $B_2$ together
(or vice-versa, depending on the position  of largest element in $\{A\vbar B\}$).
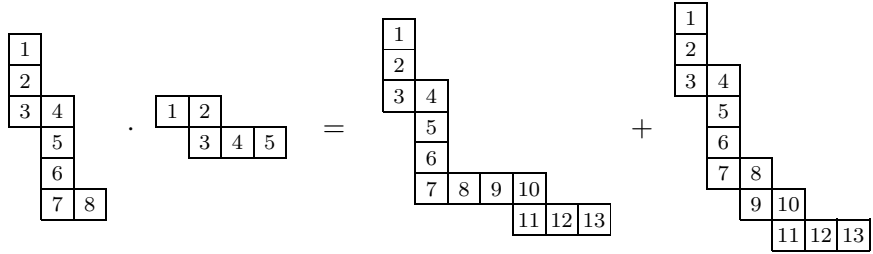
  $$\begin{matrix} 
               \Young{\Box{1} \cr \Box{2}\cr  \Box{3} & \Box{4}\cr
                \cr & \Box{5} \cr & \Box{6}   \cr  &\Box{7} & \Box{8} \cr } 
     \end{matrix}    
      \,\,\, \cdot  \,\,\,
     \begin{matrix} 
               \Young{\Box{1}  & \Box{2} \cr & \Box{3} & \Box{4} &
                \Box{5} \cr  }  
     \end{matrix}    
                \,\,\,\,\,\,=\,\,\,\,\,\,
       \begin{matrix} 
               \Young{\Box{1} \cr \Box{2}\cr  \Box{3} & \Box{4}\cr
                \cr & \Box{5} \cr & \Box{6}   \cr  &\Box{7} & \Box{8} &
                \Box{9}  & \Box{10} \cr  & & & & \Box{11} & \Box{12} &
                \Box{13} \cr   } 
     \end{matrix}      
           \,\,\, + \,\,\,
      \begin{matrix} 
               \Young{\Box{1} \cr \Box{2}\cr  \Box{3} & \Box{4}\cr
                \cr & \Box{5} \cr & \Box{6}   \cr  &\Box{7} & \Box{8} \cr
              & &   \Box{9}  & \Box{10} \cr  & & & \Box{11} & \Box{12} &
                \Box{13} \cr  } 
     \end{matrix}   
     $$ 
\end{proof}

%

This algebra isomorphism also arises as a quotient space.  Define the
two sided ideal of $\NCSI$
 generated by the monomials $\{ \M_A \,|\, \ell(A) \geq 3\}$ as $I_3$.  Notice
 that if $A$ has $\ell(A) \geq 3$ then every term in the product $\M_A \M_B$
 will be indexed by a set partition of length greater than or equal to $3$
 and hence the ideal is linearly spanned by this set of monomials as well.
 We have then that $\NCSI \slash I_3 \simeq \NSym$ since the quotient will
 be linearly spanned by the $\M_A$ for $\ell(A) \leq 2$.

There is another closely related copy of $\NSym$ sitting inside $\NCSI$.
Let $\alpha$ be a composition, and let $A(\alpha)$ be the corresponding
 canonical set partition from equation \eqref{eq:compsp}.
 We define $M_{\alpha}$ to be the sum of all 
monomials in $\NCSI$ indexed by those set partitions $A$ that can be
 obtained from $A(\alpha)$ by gluing nonconsecutive blocks. For instance, 
if $\alpha= ( 2, 1, 3, 2 )$, then
$A(\alpha)=\{1 2 \vbar 3 \vbar 4 5 6 \vbar7 8 \}$, and we can only obtain the following 
five set partitions
$\{ 1 2 \vbar 3 \vbar 4 5 6 \vbar 7 8  \}$, $ \{ 1 2 4 5 6 \vbar 3 \vbar 7 8  \}$, 
$ \{ 1 2 \vbar 3 7 8 \vbar 4 5 6  \}$, $\{ 1 2 4 5 6 \vbar 3 7 8  \}$, $\{1 2 7 8\vbar 3\vbar 4 5 6\}$. A second 
way of describing $M_{\alpha}$ is as the sum of all monomials in 
$\bbQ \langle x_1, x_2, \ldots \rangle$ whose exponents are given by composition $\alpha$. 
For instance, $M_{(2, 1, 3, 2) }= \sum_{i \neq j \neq k \neq l} x_i^2x_jx_k^3x_l^2$
 where in the sum we allow any of the possibilities of $i=k$, $i=\ell$ or $j=\ell$.

\begin{prop} 
The map $\zeta: \NSym \rightarrow \NCSI$ by
 $\zeta(R_\alpha) = M_\alpha$ is an injective algebra homomorphism.
 \end{prop}
 
 \begin{proof} The map
 $\zeta$ is clearly injective, hence it
 suffices to show that $\zeta( \R_\alpha \R_\beta ) = \zeta( \R_\alpha) \zeta( \R_\beta)$.
 We have that
 \begin{align*}
 \zeta( \R_\alpha \R_\beta ) &= \zeta( \R_{\alpha\rhd\beta} + \R_{\alpha \cdot \beta})\\
  &= \sum_{i_1 \neq   \cdots \neq i_{\ell(\alpha)} = j_1 
  \neq \cdots \neq j_{\ell(\beta)}} x_{i_1}^{\alpha_1}
 \cdots x_{i_\ell(\alpha)}^{\alpha_{\ell(\alpha)}} x_{j_1}^{\beta_1}
 \cdots x_{j_\ell(\beta)}^{\beta_{\ell(\beta)}}\\
 &+\sum_{i_1 \neq   \cdots \neq i_{\ell(\alpha)}
  \neq j_1 \neq \cdots \neq j_{\ell(\beta)}} x_{i_1}^{\alpha_1} 
 \cdots x_{i_\ell(\alpha)}^{\alpha_{\ell(\alpha)}} x_{j_1}^{\beta_1} 
 \cdots x_{j_\ell(\beta)}^{\beta_{\ell(\beta)}}\\
 &= \left( \sum_{i_1   \neq \cdots \neq i_{\ell(\alpha)}} x_{i_1}^{\alpha_1}  
 \cdots x_{i_\ell(\alpha)}^{\alpha_{\ell(\alpha)}} \right)
 \left( \sum_{j_1   \neq \cdots \neq j_{\ell(\beta)}} x_{j_1}^{\beta_1}  
 \cdots x_{j_\ell(\beta)}^{\beta_{\ell(\beta)}} \right)\\
 &=\zeta(\R_\alpha) \zeta(\R_\beta).
 \end{align*}
 \end{proof}

The following observation, due to Florent Hivert \cite{H}, shows us that $\NSym$ is also a quotient of $ \NCSI$. 

When $A$ is not equal to $A(\alpha)$ for any $\alpha$ then we will say 
that $A$ has crossings.  We remark that if $A$ has crossings then
so will every term in the expansion of $\M_{A} \M_{B}$.  Consider the two
sided ideal $\cideal$ generated by all $\M_{A}$ such that $A$ has crossings.
This ideal is then linearly spanned by all $\M_{A}$ such that $A$ has
crossings.

Now consider the quotient $\NCSI/\cideal$.  It is linearly spanned by
the basis $\M_{A(\alpha)}$ for $\alpha$ a composition.  The proof that
the elements $\NSaNC( \h_\alpha )$ are all linearly independent also
shows that they will be linearly independent in the quotient $\NCSI/\cideal$.

\begin{corollary}
\begin{equation*}
\NSym \simeq \NCSI/\cideal.
\end{equation*}
as algebras.  The isomorphism is given explicitly as $\rho : \NCSI \rightarrow
\NSym$ by 
\begin{equation*}
\rho( \M_{A} ) = \left\{ \begin{array}{cl}
\M_{A} & \hbox{ if $A$ has no crossings }\\
0 & \hbox{ otherwise }
\end{array}\right..
\end{equation*}
\end{corollary}

A computation of $\Delta^{\NCSI} \circ \rho \circ \NSaNC( \h_{3} )$
and $(\rho \otimes \rho) \circ \Delta^{\NCSI} \circ \NSaNC( \h_{3} )$ 
shows that these spaces are not isomorphic as Hopf algebras since $I$ is
 not a Hopf ideal.

\end{subsection}

\end{section}
 
 \begin{section}{The Harmonics with respect to the Hausdorff derivative.} \label{sec:lie}

We give an elegant characterization of  the space of harmonics in noncommuting 
variables with respect to the Hausdorff derivative in terms of the free Lie algebra.
We will require some basic definitions and results for which we refer the reader 
to \cite{R} for references and their proofs.

A Lie algebra over $\C$ is a $\C$-module $\L$, together
with a bilinear mapping
\begin{align*}
\L \times \L &\to \L\\
(x,y)  &\mapsto [x,y]
\end{align*}
called the Lie bracket.  This bracket must satisfy two
identities,
$[x,y] = -[y,x]$
and
$[x,[y,z]]+[y,[z,x]] + [z,[x,y]] = 0.$
Subalgebras of Lie algebras, homomorphisms and modules are
defined as usual for Lie algebras. 
Any associative algebra $\A$ over $\C$ acquires a natural structure of a
Lie algebra when $[x, y]$ is defined by $[x,y]=xy-yx$. 
The free Lie algebra
can be realized as the linear span of the
minimal set of polynomials in $\polys$
which include the variables $\{ x_1, x_2, \ldots, x_n \}$ and
is closed under the bracket operation.

For a Lie algebra $\L\subseteq \polys$ with the natural bracket operation, 
the enveloping algebra of $\L$
is the subalgebra of $\polys$ generated by 
the elements of $\L$ under the concatenation product.

Let $\L=\L(X_n)$  be the free Lie algebra generated by the noncommutative
alphabet $X_n=\{x_1, x_2, \cdots, x_n\}$, and let $\L'=[\L, \L]$ be the Lie subalgebra generated by the brackets $[P,Q]$ where both $P$ and $Q$ are in $\L$. Let $\A'$ be the enveloping  algebra of $\L'$. In particular, $\L= \L' + \C X_n$,
where $\C X_n$ denotes the linear polynomials.

We want to characterize the harmonics of the symmetric group in noncommuting
variables. Recall that in the commutative setting the harmonics are defined
as the set of solutions for the system of PDE obtained by looking at
symmetric functions as differential operators. 

Our goal is to compute  the harmonics of the symmetric group in the 
noncommutative setting. To this end, we should  start by defining what 
we mean by the derivative of a noncommutative polynomial. We  first focus
 our attention on the Hausdorff derivative, the most common definition for derivative in
the noncommutative setting.

Let $w$ be a monomial in $\C \langle X \rangle$, that is, a word.
The Hausdorff derivative of $w$ with regard to the letter $x$ is defined as
the sum of all  subwords $w'$  obtained from $w$  by  deleting an occurrence of letter $x$,
and then extended by linearity. 
For instance,
$ \haus_x xyx^2y= yx^2y + 2 xyxy$, and $\haus_x [x,y]= \haus_x (xy) - \haus_x (yx)=0$.
 
The following theorem can be found in \cite{R} and 
characterizes the elements of $\A'$ as the elements of $\polys$ that are killed by each derivation.

\begin{prop} (\cite{R})
\[
\bigcap_{x\in X_n} \ker{\haus_x} = \A'~.
\]
\end{prop}
%


For any polynomial $f \in \polys$, we will denote by
$f(\haus_{X_n})$ the linear differential
operator formed by replacing each of the monomials $x_{i_1} x_{i_2} \cdots x_{i_k}$
by the differential $\haus_{x_{i_1}} \haus_{x_{i_2}} \cdots \haus_{x_{i_k}}$.
Note that $\haus_x \haus_y = \haus_y \haus_x$ and so we have that 
the operator $\M_A(\haus_{X_n})$ acts up to constant as $m_{\la(A)}(\haus_{X_n})$.
More precisely, $\M_A(\haus_{X_n})( f(X_n) ) = \chi(\M_A)(\haus)(f(X_n))$.

\begin{defn}  Let $X_n = \{x_1, x_2, \ldots, x_n\}$ be a finite noncommuting alphabet,
the harmonics  with respect to the Hausdorff 
derivative are defined  as the space of solutions of the system of PDEs 
$$f(\haus_{X_n}) Q(X_n) = 0$$
for all $f \in \NCSI_n$ without constant term.  We denote the solution space by $\MH_n.$
\end{defn}

%

\begin{theorem}  (Poincar\'e-Birkhoff-Witt) Let $\L$ be a Lie algebra and
consider  $\L$ as a vector space  with a totally ordered basis $(w_i)_{i \in I}$.
Let $\A_0$ be its enveloping algebra and $\varphi_0: \L \to \A_0$ be the natural
Lie algebra homomorphism. Then $\A_0$ is a vector space over $\C$ with basis
$\varphi_0(w_{i_1}) \ldots \varphi_0(w_{i_n}),$ where $n \ge 0$, $i_1, \ldots, i_n \in I$,
and $i_1\ge \ldots \ge i_n.$
\end{theorem}

Let $\L$ be the free Lie algebra in the variables $X_n$.
Take a basis  $\mathcal{B}'$ of $\L' = [\L, \L]$,  then
the linear polynomials $\C X_n$ satisfy $\L = \L' \oplus \C X_n$ and 
\[
\mathcal{B}=\mathcal{B}' \cup X_n
\]
is a basis for $\L$.
Next, order the basis in such a way that the elements of $\mathcal{B}'$ 
are strictly bigger than the elements of $X_n$.

The enveloping algebra of the free Lie algebra $\L$ is $\polys$. 
Therefore, the theorem of Poincar\'e-Birkhoff-Witt  implies that 
decreasing products of elements of $\mathcal{B}$ form a basis of $\polys$.
Moreover, since $\mathcal{A'}$ is the enveloping algebra of $\L'$, the 
theorem of Poincar\'e-Birkhoff-Witt also implies that decreasing products 
of $\mathcal{B}'$ are a basis of $\A'$.  We also know that decreasing products of
$X_n$ are isomorphic to $\C[X_n]$. We conclude that, as vector spaces,
\begin{equation} \label{iso}
\polys \simeq  \A' \otimes \C[X_n] \simeq \A'[X_n]
\end{equation}
Furthermore, this isomorphism is compatible with derivations $\haus_a$.
That is, for a $P(X_n) \in \polys$ where $P(X_n) = \sum_{i} b_i f_i(X_n)$ 
and $b_i \in \A'$ and $f_i(X_n) \in
\C[X_n]$, we have
$$\haus_{a} P(X_n) = \sum_i b_i \haus_{a} f_i(X_n)$$
for $a \in X_n$.  This follows because $\haus_a (\A')=0$.

We have from this discussion the following theorem.
\begin{theorem} \label{th:mharmiso} 

Let ${\mathcal H}_n$ be the classical harmonics. That is,
$${\mathcal H}_n = \{ f(X_n) \in \C[X_n] : p( \haus_{X_n} ) f(X_n) = 0 \hbox{ for all } 
p( X_n ) \in \Sym_n \hbox{ with } p(0) = 0\}.$$

Then, as vector spaces,
$$\MH_n   \simeq \A'  \otimes {\mathcal H}_n.$$

\end{theorem}

Moreover, Chevalley \cite{C} showed that 
${\mathcal H}_n$ the linear span of derivatives of the Vandermonde polynomial 
$$\Delta_n = \prod_{1 \leq i < j \leq n} (x_i - x_j) 
= \sum_{\pi \in \s_n} sgn(\pi) \,\, x_{n}^{\pi_1-1} x_{n-1}^{\pi_2-1} \cdots x_{1}^{\pi_n-1}.$$

Note that  from the existence of the isomorphism \eqref{iso}
and the classical characterization for the harmonics in the 
commutative case, we obtain that
$N\Delta_n \in \MH_n$, where $N\Delta_n$ is the noncommutative Vandermonde, defined as
\[
N\Delta_n = \sum_{\pi \in \s_n} sgn(\pi) \,\, x_{n}^{\pi_1-1} x_{n-1}^{\pi_2-1} \cdots x_{1}^{\pi_n-1}.
\]
It is interesting to note that 
all other possible noncommutative Vandermondes (obtained by fixing an 
order in the variables) are also harmonics, but one 
suffices to describe this space. To see this, we only need
to order $X_n$ in all possible ways 
before applying the  Poincar\'e-Birkhoff-Witt  theorem.  Likewise, we have that
each derivative of $N\Delta_n$ is in $\MH_n$. Therefore, 
$Span_{\haus}[N\Delta_n] \subseteq \MH_n$.  Theorem \ref{th:mharmiso} implies
that $\MH_n$ is equal to the $\A'$-module generated by $N\Delta_n$ and all
its derivatives.

%


\label{sec:chevalley}
A famous theorem due to Chevalley says that the ring of polynomials is
isomorphic to the tensor product of its invariants times its coinvariants (that 
in the commutative case are show to be isomorphic to the harmonics). 
\begin{theorem}[Chevalley \cite{C}] \label{th:chev} As $\s_n$-modules,
\[
{\mathcal H}_n \otimes \Sym_n  \simeq \C[X_n]~.
\]
\end{theorem}

We conclude a mixed commutative/non\-com\-mu\-ta\-tive 
version of Che\-val\-ley's theorem which holds on the level
of vector spaces which we derive from the results above. But to get the isomorphism as $\s_n$-module
we need some more tools. 

We first define on $\polys$ a commutative product. 
The shuffle product, denoted by $\shuf$, is the bilinear operation recursively defined as follow.
Given variables $x,y$ and monomials $v, w \in \polys$,
  $$1\shuf u = u\shuf 1 =u \quad\hbox{and}\quad xu\shuf yv=x(u\shuf yv) + y(xu\shuf v).$$
This is a well known commutative product on $\polys$.
It is clear that the forgetful map $\chi: \polys \rightarrow {\mathbb Q}[X_n]$ acts as
  $$\chi(x_{i_1}\shuf x_{i_2}\shuf \cdots \shuf x_{i_k})= k! x_{i_1} x_{i_2}  \cdots  x_{i_k}.$$
Define now $\tilde{p}_k=\sum_{i=1}^n x_{i}\shuf x_{i}\shuf \cdots \shuf x_{i} = \M_{\{[k]\}}[X_n]$
where the variable $x_i$ is shuffled with itself $k$ times. For $\lambda=(\lambda_1,\lambda_2,\ldots,\lambda_{\ell(\la)})\vdash m$ a partition of the integer $m$ we let $\tilde{p}_\lambda=\tilde{p}_{\lambda_1}\shuf\cdots\shuf \tilde{p}_{\lambda_{\ell(\la)}}$. We then have that
\begin{equation}\label{chiP}
  \chi(\tilde{p}_\lambda)= \la! {p}_{\lambda_1}{p}_{\lambda_2}\cdots {p}_{\lambda_{\ell(\la)}},
\end{equation}
where $p_k=\sum_{i=1}^n x_i^k\in\Sym_n$ is the classical power sum symmetric polynomial. If we denote by $\widetilde{\Sym}_n\subseteq \polys$ the vector space spanned by the $\tilde{p}_\lambda$
with $1\le \lambda_i\le n$.
\begin{lemma} \label{TSS} As graded $\s_n$-modules, 
  $ \widetilde{\Sym}_n \simeq \Sym_n$. 
\end{lemma}

\begin{proof} It is well known that $\Sym_n=\bbQ[p_1,p_2,\ldots,p_n]$. 
The map $\chi$ in equation~(\ref{chiP}) restricted to $\widetilde{\Sym}_n$ gives
us a surjective linear map $\chi\colon \widetilde{\Sym}_n\to \Sym_n$. This 
map preserves the degree of the polynomial, so we can restrict our attention
to the homogeneous component of degree $m$, $\widetilde{\Sym}_n^{m}$.
Since the product $\shuf$
is commutative, $\dim(\widetilde{\Sym}_n^{m})\le \dim(\Sym_n^{m})$ the 
number of partitions $\lambda\vdash m$ with $1\le \lambda_i\le n$.
Hence $\chi : \widetilde{\Sym}_n \rightarrow \Sym_n$ is an isomorphism 
of graded vector spaces. Since each element of $\widetilde{\Sym}_n$ is
$\s_n$ invariant (as is $\Sym_n$), $\widetilde{\Sym}_n$ and $\Sym_n$ are
isomorphic as $\s_n$ modules as well.
\end{proof}

Let us denote by $\NCSI_n^+$ the set of $f\in\NCSI_n$ without constant term.
Recall that the map $\chi\colon\NCSI_n\to\Sym_n$ is surjective and also that the Hausdorff derivative commutes. We thus have $f(\partial)=\chi(f)(\partial)$ for all $f\in\polys$. Combining these remarks, we get
 \begin{align*}
   \MH_n &=\{P\in\polys \,|\, f(\partial)P=0,\  \forall f\in\NCSI_n^+\} \\
               &=\{P\in\polys \,|\,\chi( f)(\partial)P=0, \ \forall f\in\NCSI_n^+\}\\
               &=\{P\in\polys \,|\, p_k(\partial)P=0, \ 1\le k\le n\}.
\end{align*}
Now let $\langle\, , \rangle$ denote the scalar product on $\polys$ for which the monomials forms an orthonormal basis. For all variable $x\in X_n$ and monomials $u,v\in\polys$ we easily see that
\begin{equation} \label{scshuf}
  \langle x\shuf u,v\rangle = \langle u,\partial_x v\rangle.
\end{equation}
Finally let $\langle \tilde{p}_k : {1\le k\le n}\rangle_\shuf\subseteq\polys$ denote the ideal generated using the shuffle product. That is
  $$ \langle \tilde{p}_k : {1\le k\le n}\rangle_\shuf = \left\{ \sum_{k=1}^n \tilde{p}_k\shuf q_k \,|\, q_k\in \polys\right\}.$$

\begin{lemma} $$\MH_n=\langle \tilde{p}_k : {1\le k\le n}\rangle_\shuf^\perp.$$.\end{lemma}

\begin{proof} If $P\in \MH_n$ then for all $1\le k\le n$ we have $p_k(\partial)P=0$.                                                                                         Given any $F\in \langle \tilde{p}_k : {1\le k\le n}\rangle_\shuf$, 
$F = \sum_{k=1}^n p_k \shuf q_k$ where
$q_k \in polys$ and we calculate
  \begin{align*}
   \langle F, P\rangle  &=\sum_{k=1}^n \sum_{i=1}^n \langle x_i\shuf \cdots \shuf x_i\shuf q_k, P\rangle\\
               &=\sum_{k=1}^n \sum_{i=1}^n \langle  q_k, \partial_{x_i}^k P\rangle 
               =\sum_{k=1}^n  \langle  q_k, p_k(\partial) P\rangle = 0,
\end{align*}
where we have use the identity~(\ref{scshuf})  $k$ times in the summands. Hence $P\in \MH_n$ implies $P\in \langle \tilde{p}_k : {1\le k\le n}\rangle_\shuf^\perp$. Conversely if  
$P\not\in \MH_n$, then there is a $1\le k\le n$ such that $p_k(\partial)P\not = 0$. This means we can find $q\in\polys$ such that
  $$0\not = \langle q,p_k(\partial)P\rangle =\sum_{i=1}^n \langle q,\partial_{x_i}^k P\rangle
    = \langle \tilde{p}_k\shuf q,P\rangle$$
and conclude that $P\not\in  \langle \tilde{p}_k : {1\le k\le n}\rangle_\shuf^\perp$.
\end{proof}

At this point, we have shown that $\polys = \MH_n\oplus\  \langle \tilde{p}_k : {1\le k\le n}\rangle_\shuf$. This gives us for any $G\in\polys$
\begin{equation}\label{recG}
  G=P + \sum_{k=1}^n \tilde{p}_k\shuf q_k = 1\shuf P + \sum_{k=1}^n \tilde{p}_k\shuf q_k,
\end{equation}
where $P\in \MH_n$ and $\deg(q_k) <\deg(G)$. If we repeat the use of equation~(\ref{recG}) recursively on the $q_k$ we get that
  $$G=\sum_\lambda \tilde{p}_\lambda \shuf P_\lambda,$$
where the sum runs over $\lambda=(\lambda_1,\ldots,\lambda_{\ell(\la)})$ such that $1\le \lambda_i\le n$ and $P_\lambda\in\MH_n$. Also by convention we allow $\lambda=()$ and $\tilde{p}_{()}=1$.
This equation shows that the graded linear map
\begin{equation}\label{Phimap}
\psi\colon \widetilde{\Sym}_n\otimes \MH_n \to \polys,
\end{equation}
defined by $\psi(\tilde{p}_\lambda\otimes P)= \tilde{p}_\lambda\shuf P$, is surjective.

\begin{theorem} As graded $\s_n$-modules,
\[
\Sym_n \otimes \MH_n
\simeq \polys~.
\]
\end{theorem}
\begin{proof} 
By equation \eqref{iso} and Theorems \ref{th:mharmiso} and
\ref{th:chev}, we have
\begin{align*}
 \Sym_n \otimes  \MH_n \simeq \Sym_n \otimes {\mathcal H}_n  \otimes  \A' 
\simeq \C[X_n]  \otimes  \A' 
\simeq \polys.
\end{align*}
as vector spaces. Combined with Lemma~\ref{TSS} and the surjectivity of $\psi$,  this shows that 
  $$\psi\circ(\chi^{-1}\otimes id) \colon \Sym_n \otimes  \MH_n \to \widetilde{\Sym}_n \otimes \MH_n \to
 \polys$$
is surjective, and hence an isomorphism of vector spaces.  

To view the result as an isomorphism of $\s_n$-modules, we first need to make 
sure that $\MH_n$ is indeed an $\s_n$-module. This follows from the fact that for 
any $P\in\MH_n$, any $\sigma\in\s_n$ and for all $1\le k\le n$, 
$p_k(\partial) \sigma(P)=\sigma (p_k(\partial) P)=0$ and thus $\sigma(P)\in\MH_n$.
We already know that $\chi$ restricted to $\widetilde{\Sym}_n$ is a morphism of $\s_n$-modules. 
It thus remains to show that $\psi \circ (\chi^{-1} \otimes id)$ is also a morphism of $\s_n$-modules. 
For this let $\sigma\in\s_n$:
$$ \sigma\circ\psi(\tilde{p}_\lambda\otimes P)=\sigma( \tilde{p}_\lambda\shuf P)=(\sigma\tilde{p}_\lambda)\shuf (\sigma P)=\psi((\sigma\tilde{p}_\lambda)\otimes(\sigma P))
 = \psi\circ(\sigma \otimes \sigma) (\tilde{p}_\lambda\otimes P).
$$
Since both $\tilde{p}_\lambda$ and $p_\la$ are both $\s_n$ invariant,
$\sigma\circ\psi \circ (\chi^{-1} \otimes id)= \psi \circ (\chi^{-1} \otimes id) \circ (\sigma \otimes \sigma)$
and our proof is then complete.
\end{proof}

\end{section}

 \begin{section}{Noncommutative invariants of the symmetric group} \label{sec:Wolfsummary}

In classical invariant theory of the symmetric group (see \cite{Mac,
steinberg}) the ring of symmetric polynomials in $n$ (commuting) variables
is free. In particular it is a polynomial ring with $n$ generators, one in
each degree. In \cite{W}, Wolf was the first to study $\NCSI_n$ as
invariants in noncommuting variables. Her main theorem shows that the space
of noncommutative invariants of the symmetric group is also free. 

In her proof, it is not obvious how to construct the generators and the
combinatorics of set partitions is not fully developed. In particular, it is
not clear what the Hilbert series of the invariant polynomial ring in
$n$ noncommutative variables is.
In the final section of this article we will need her result and the associated Hilbert
series. We thus present it here along with a constructive proof.

Given two set partitions $A=\{A_1, \ldots , A_k\} \vdash [n]$ and
$B=\{B_1,\ldots , B_{\ell}\} \vdash [m]$, we define
\begin{equation*}
A\circ B= \left\{ \begin{array}{cl}
\{A_1\cup (B_1+n), \ldots, A_k\cup (B_k+n), (B_{k+1}+n), \ldots,
(B_{\ell}+n)\} & \hbox{ if $k\le \ell$ }\\ \\
\{A_1\cup (B_1+n), \ldots, A_{\ell}\cup (B_{\ell}+n), A_{{\ell}+1}, \ldots,
A_k\} \hfill& \hbox{ if $k>\ell$ }
\end{array}\right..
\end{equation*}
Recall that the parts of $A$ and $B$ are ordered according to the minimum
elements in each part.
For example, if  $A=\{13\vbar2\}$ and $B=\{1\vbar2\vbar3\}$ then $A\circ B=\{134\vbar25\vbar6\}$
and $B\circ A=\{146\vbar25\vbar3\}$.
We note that $\ell(A\circ B) = \max({\ell(A),\ell(B)})$.

If $A=B\circ C$ for $B$ and $C$ nonempty set partitions, 
then we say that $A$ {\sl splits}. If it is not possible to
split $A$, then we say that it is {\sl nonsplitable}.  By convention,
only non-empty set partitions are nonsplitable.

\begin{exem}  For $n=3$. the list of all set partitions is
\begin{align*} 
 &\{123\} = \{ 1 \}\circ  \{ 1 \} \circ  \{ 1 \}
 &&\{1\vbar23\} \hbox{ nonsplitable }
 &&\{13\vbar2\} = \{1\vbar2\}\circ \{1\}\\
 &\{12\vbar3\} = \{ 1 \}\circ  \{ 1\vbar2 \}
 &&\{1\vbar2\vbar3\} \hbox{ nonsplitable }
 &&\\
 \end{align*}
 \end{exem}

As we remarked in Section \ref{sec:Hopfalg} a basis for $\NCSI_n$ is given by
$\{\M_A[X_n]\}_{\ell(A)\le n}$.
Consider the set $X_n$ as an alphabet where $x_1<x_2<\cdots<x_n$. A monomial
in these noncommutative variables can be viewed as a word in the alphabet
$X_n$.
Given $A$ such that $\ell(A)\le n$, we order the monomials of $\M_A[X_n]$ by
lexicographic order and denote by $LT(\M_A[X_n])$ the smallest
monomial in $\M_A[X_n]$. For example,
$LT(\M_{\{14\vbar25\vbar3\}}[X_6])=x_1x_2x_3x_1x_2$. In general, the $k$th variable
of $LT(\M_A[X_n])$ is $x_i$ exactly when $k\in A_i$ and the parts of $A$
are ordered according to the minimum elements in each part.

\begin{lemma}  For any set partitions $A$ and $B$ with at most $n$ parts, we
have
\begin{equation}\label{leadterm}
 LT\big(\M_A[X_n]\M_B[X_n]\big) =  LT(\M_A[X_n]) LT(\M_B[X_n])~.
\end{equation}
\end{lemma}

\begin{proof}
This is a direct consequence of the following well known fact about
lexicographic order. Given four words (monomials) $u_1,u_2,v_1,v_2$ in the
alphabet $X_n$ such that $u_1\le_{lex} v_1$ and $u_2\le_{lex} v_2$, then
$u_1u_2\le_{lex}v_1v_2$. So the smallest term of $\M_A[X_n]\M_B[X_n]$
is the product of the smallest term of $\M_A[X_n]$ with the smallest term of
$\M_B[X_n]$, and all terms have coefficient equal to $1$.
\end{proof}

We now proceed to show the main results of Wolf.

\begin{prop} $\NCSI_n$ is freely generated as an algebra by
  $$\{\M_A[X_n]\ :\ \ell(A)\le n \hbox{ and $A$ is nonsplitable}\}.$$
\end{prop}

\begin{proof}
Let $B\vdash [m]$ be a  set partition such that $\ell=\ell(B)\le n$. Let
$1 \le k\le m$ be the smallest integer such that
$$min( B_i \cap [k]^c) \leq min(B_j \cap [k]^c)\hbox{ for each $1\leq i<j\leq \ell(B)$}$$
where we use the convention that $[k]^c = \{ k+1, \ldots, m\}$ and $min( \emptyset )  = \infty$.

Let $B^{(1)}=\{B_i\cap [k]\,:\, i\le r\}\vdash [k]$, then by the choice of $k$,
$B^{(1)}$ is nonsplitable. 

If $k=m$, then $B=B^{(1)}$ is nonsplitable, otherwise
$k<m$ and $\tilde{B}=\{(B_i\cap\{k+1,\ldots,m\})-k\,:\, i\le
\ell\}\vdash[m-k]$ with  $B=B^{(1)}\circ \tilde{B}.$
 Repeating this process recursively, we obtain a {\sl unique} decomposition
of $B$ into nonsplitable set partitions:
  $$ B=B^{(1)}\circ B^{(2)}\circ \cdots \circ B^{(s)}. $$
 Now consider the following expansion:
 \begin{equation}\label{eq:Wbasis}
 W_B[X_n] := \M_{B^{(1)}}[X_n]\M_{B^{(2)}}[X_n]\cdots \M_{B^{(s)}}[X_n]=\sum_{D\vdash
[m]} c_{D} \M_D[X_n].
\end{equation}

Since $LT(\M_{B^{(1)}}[X_n]\cdots \M_{B^{(s)}}[X_n])=LT(\M_B[X_n])$ we  have
that $c_B=1$ and $c_D=0$ whenever $LT(\M_D[X_n])<_{lex}LT(\M_B[X_n])$.
This implies that the change of basis matrix between the basis
$\{\M_B[X_n]\,:\, \ell(B)\le n\}$ and the set ${\mathcal W} = \{ W_B[X_n] \,:\, \ell(B) \le n\}$ is
upper triangular and therefore ${\mathcal W}$ is a 
linear basis of $\NCSI_n$. 
We now remark that $\mathcal W$ is
also a basis of the free noncommutative algebra generated by the set
$\{\M_A[X_n]\ :\ \ell(A)\le n \hbox{ and $A$ is nonsplitable}\}$ and this
concludes the proof.
\end{proof}

Let $S_{m,k}$ denote the number of set partitions of $m$ with exactly $k$ parts
(the Stirling numbers of the second kind).
Then the number of set partition of $m$ with at most $n$ parts is
$\sum_{k=1}^n S_{m,k}$. We thus have that
   $$B_n(q)=\sum_{m\ge 0} \dim_m(\NCSI_n) q^m = \sum_{m\geq0}
\sum_{i=1}^n S_{m,i} q^m,$$
 where $dim_m(\NCSI_n)$ is the dimension of the
homogeneous component of degree $m$ in $\NCSI_n$.
 Let $w_{m,n}$ be the number of nonsplitable set partitions of $m$ with at most
$n$ parts and let $W_n(q)=\sum_{m\ge 0} w_{m,n}q^m$. A direct consequence of
the previous theorem is that
 $B_n(q)= (1-W_n(q))^{-1}$. Thus
 \begin{equation}\label{Hnonsplit}
 W_n(q) = 1 - \frac{1}{B_n(q)}.
 \end{equation}

As we will require the use of these numbers later, we include a table of
the values of $w_{m,n}$ for $1 \leq m,n \leq 8$.

\medskip

\begin{equation*}
\begin {array}{cccccccccc}
m\slash n&\vline&1&2&3&4&5&6&7&8\\
\hline
1&\vline&1&1&1&1&1&1&1&1\\
2&\vline&0&1&1&1&1&1&1&1\\
3&\vline&0&1&2&2&2&2&2&2\\
4&\vline&0&1&5&6&6&6&6&6\\
5&\vline&0&1&13&21&22&22&22&22\\
6&\vline&0&1&34&78&91&92&92&92\\
7&\vline&0&1&89&297&406&425&426&426\\
8&\vline&0&1&233&1143&1896&2119&2145&2146
\end {array}
\end{equation*}

\medskip 

Notice that the differences between adjacent entries in this tables
is the number of generators of a fixed length as described by Wolf.
\end{section}

\begin{section}{The coinvariants of the symmetric group}
 \label{sec:coinvariants}
 
 Let $X_n$ be an alphabet with $n$ letters, and let
 $\polys$ be the corresponding
ring of noncommutative polynomials. 
%



%

%


We denote by $\NCSIdeal = 
{\mathcal L}\left\{ P(X_n) \M_A[X_n]  \,|\, k \geq 1, A \vdash [k], P(X_n) \in \polys \right\}$
the left ideal of $\polys$ generated by all
elements of $\NCSI$ without constant term.
The coinvariant algebra of the
symmetric group in noncommutative 
variables will be defined as the quotient:
\[
  \polys / \NCSIdeal~.
\]

To find a linear basis of the space $\polys/\NCSIdeal$ we
use some standard techniques of the
theory of languages, that can be found in \cite{BR}.
We start by introducing some definitions.

%


Let $L^*$ be the free monoid generated by $L$, an alphabet for the monoid.
A {\sl suffix set} is a subset $C$ of $L^*$ such that 
for all $u$ and $v$ in $L^*$, if $v$, $uv$ are both in
$C$ implies  that $u=\emptyset$.  

A subset $P$ of $L^*$ is {\sl prefix closed} (resp. {\sl suffix closed}) if $uv \in P$
(resp. $vu \in P$) implies $u \in P$ for all words $u$ and $v$.
There is a bijection between suffix sets and suffix
closed sets. To a suffix set $C$ is associated the
suffix-closed set $P=L^* \setminus L^* C$, that is the
 set of words which do not end with an element of $C$.
Moreover,  $L^* = P C^*$.

For a polynomial $P(X_n) \in \polys$, the leading term with
respect to the lexicographic order will be denoted $LT(P(X_n))$
(without the leading coefficient, hence
$LT(P(X_n))$ will be an element of $X_n^\ast$).

Noncommutative monomial symmetric functions 
indexed by set partitions have the property that
$\{ LT(\M_A[X_n]) \,|\, A \hbox{ a set partition}\}$ is a
prefix closed set.  That is, any prefix $u$ of $LT( \M_A[X_n])$
is the leading term of $\M_B[X_n]$ for some set partition $B$.
In particular $B$ will be equal to $A$ restricted to
$\{ 1,2, \ldots, |u|\}$.
For example, if $A = \{ 13\vbar 246\vbar 5\}$ has leading term
$x_1 x_2 x_1 x_2 x_3 x_2$ and the heads of these monomials
are $x_1$, $x_1 x_2$,
$x_1 x_2 x_1$, $x_1 x_2 x_1 x_2$, 
$x_1 x_2 x_1 x_2 x_3$ corresponding to the leading
terms of $\M_{\{1\}}[X_n]$, $\M_{\{1\vbar2\}}[X_n]$,
$\M_{\{13\vbar2\}}[X_n]$,$\M_{\{13\vbar24\}}[X_n]$, and $\M_{\{13\vbar24\vbar5\}}[X_n]$.

We established in the previous section
that the $W_A[X_n]$ of equation \eqref{eq:Wbasis}
for $A$ a set partition with $\ell(A) \leq n$
are a basis for $\NCSI_n$ we will use this to show that 
the polynomials $u\,\M_A[X_n]$ for $A$ nonsplitable 
and $u \in X_n^\ast$ form a basis for $\NCSIdeal$.

\begin{prop} \label{nspsuffix}
Let
$$C = 
\{LT(\M_{A}[X_n]) \,|\, A\hbox{ nonsplitable}, \ell(A) \leq n\}.$$
$C$ is a suffix set of the language $X_n^\ast$.
\end{prop}

\begin{proof}
Suppose $u$ and $vu$ are both in $C$ and $v \neq \emptyset$.  Then
since $vu = LT( \M_A[X_n] )$ for $A$ nonsplitable
then $v = LT(\M_B[X_n])$ for
some nonempty set partition $B$ since it is the head of the
monomial $LT( \M_A[X_n] )$. Because $u = LT( \M_C[X_n])$
for some nonsplitable set partition $C$, we can
conclude that $A = B \circ C$, but this contradicts that $A$ is
nonsplitable.
\end{proof}

For any left ideal $I$ of $\polys$ we note that the set
$M_I = \{ LT(f) \,|\, f \in I \}$ is a left monomial ideal.  That is,
for each $v \in M_I$ and for each $u \in X_n^\ast$,
$uv \in M_I$.
%



Now we are ready to describe precisely the quotient $\polys/\NCSIdeal$.
Consider the set of leading terms of $\NCSIdeal$, $M_\NCSIdeal = 
\{ LT(f) \,|\, f \in \NCSIdeal \}$,
which is a left monomial ideal.  We note that 
\begin{align*}M_\NCSIdeal &= \{ v LT( \M_A[X_n] ) \,|\, v \in X_n^\ast,
A \hbox{ set partition}\}\\
 &= \{ v LT( \M_A[X_n] ) \,|\, v \in X_n^\ast, A \hbox{ nonsplitable set partition}\}.
\end{align*}  
That is, $M_\NCSIdeal = X_n^\ast C$ where $C$ is given
Proposition \ref{nspsuffix}.  We conclude then by the correspondence
between suffix closed sets and suffix sets that 
$$M_\NCSIdeal^c = X_n^\ast \setminus M_\NCSIdeal 
= X_n^\ast \setminus X_n^\ast C$$ 
is suffix closed.

\begin{prop} \label{prop:basisIdl}
The set
$$\{ u\,\M_A[X_n] \,|\, u \in X_n^\ast, A \hbox{ nonsplitable set partition} \}$$
is a linear basis of $\NCSIdeal$ and
any element of $w \in M_\NCSIdeal = X_n^\ast C$ can
be decomposed uniquely as $w=uv$ where $u \in X_n^\ast$ and
$v \in C$.
\end{prop}

\begin{proof}
Assume $u v = u' v' \in X_n^\ast C$
with $v, v' \in C$.
Without loss of generality assume $v = w v' \in C$. Since $C$ is a suffix set
then $w = \emptyset$ and hence $v = v'$ and $u= u'$.  Therefore the decomposition
of $u v \in M_\NCSIdeal$ is unique.

Next consider the set
$\{ u\,\M_A[X_n] \,|\, u \in X_n^\ast, A \hbox{ nonsplitable set partition} \}$.
Since the leading terms of
the elements of this set are all distinct, they are linearly independent.
This set must also span the ideal $\NCSIdeal$ because every element
of the form $v W_A[X_n] \in \NCSIdeal$ is in the linear span of this set
and the $v W_A[X_n]$ are certainly a spanning set of the ideal.
\end{proof}


\begin{prop} \label{prop:basisMIC}
 $M_\NCSIdeal^c$ is a basis for $\polys/\NCSIdeal$.
\end{prop}
\begin{proof}
First we show that $M_\NCSIdeal^c$ spans the vector space.
Since the words of $X_n^\ast$ span $\polys$, it suffices to show that for $v \in X_n^\ast$,
$$v \equiv \sum_{u \in M_\NCSIdeal^c} a_u u \mod \NCSIdeal.$$
Assume that $v$ is the smallest such monomial which is not a linear
combination of $u \in M_\NCSIdeal^c$.  Since $v \notin M_\NCSIdeal^c$, then
$v \in X_n^\ast C = M_\NCSIdeal$ and so $v = u LT( \M_A[X_n])$ for
some nonsplitable set partition $A$.
Now $v - u \,\M_A[X_n]$ is equal to a sum of terms which are smaller than 
$v$ in lexicographic order and hence are equivalent
to a linear combination of elements of $M_\NCSIdeal^c$.

The monomials in $M_\NCSIdeal^c$ are also linearly independent since if we assume that
\begin{equation} \label{MIcdesc}
P(X_n) = \sum_{u \in M_\NCSIdeal^c} a_u u \equiv 0 \mod \NCSIdeal~,
\end{equation}
then $P(X_n) \in \NCSIdeal$ and hence
$LT(P(X_n)) \in M_\NCSIdeal$.  Since the leading term of $P(X_n)$ is one of
the monomials of $M_\NCSIdeal^c$, the only way this can happen is if $P(X_n)=0$.
\end{proof}

We conclude from Propositions \ref{prop:basisIdl} and \ref{prop:basisMIC} the following corollary.

\begin{corollary}\label{cor:dimquo}
The dimension of the subspace of degree $k$ in 
$\polys / \NCSIdeal$ is equal to
\[
n^k - \sum_{i \le k}  w_{i,n}  \,n^{ k-i}
\]
where $w_{m,k}$ is the number of nonsplitable
set partitions of size $m$ and with
length less than or equal to $k$.
\end{corollary}

Define the following three generating functions:
\begin{align*}
T_n (q)&:=\sum_k dim_k (\polys)  \, q^k=\frac{1}{1-nq}\\
B_n (q)&:=\sum_k dim_k (\NCSI_n) \, q^k = \frac{1}{1-W_n(q) }=\frac{1} {1 - \sum_{k \ge 0}  w_{k,n} q^n}\\
C_n (q)&:=\sum_k dim_k (\polys/\NCSIdeal) \, q^k = \sum_{k \ge 0} \big( n^k - \sum_{i \le k} w_{i,n} n^{k-i} \big) q^k
\end{align*}

Observe that we have the relationship $T_n(q) = B_n(q) C_n(q)$ by the following
calculation.

\begin{align*}
B_n (q) C_n (q) &= B_n (q) \Big( \sum_{n \ge 0} n^k q^k \Big)
   - B_n (q) \Big( \sum_{k \ge 0} n^k q^k\Big) \Big( \sum_{k \ge 0}  w_{k,n} q^k\Big) \\
                             &= B_n(q) \Big( \sum_{n \ge 0} n^k q^k\Big) \Big( 1 - \sum_{k \ge 0}  w_{k,n} q^k \Big)\\
                             &=\sum_{n \ge 0} n^k q^k = \frac{1} {1-nq} = T_n(q).
\end{align*}

We conclude that as graded vector spaces
\[
\polys/\NCSIdeal \otimes \NCSI_n  \simeq \polys~.
\]


As mentioned in the introduction, the twisted derivative
provides a second definition of derivation in $\polys$ 
(see for example \cite{R}). It is defined as 
$$\twist_a(w) = \left\{ \begin{array}{cl} w' &\hbox{ if } w = a w' \\
0 & \hbox{ otherwise }\end{array} \right.~.$$
We can show that space of noncommutative coinvariants of the
symmetric group is isomorphic the space of harmonics
of the symmetric group with respect to the twisted derivative.

Recall that the scalar product is defined with the monomials as an
orthonormal basis.  The twisted derivative has the property analogous
to equation \eqref{scshuf} for the Hausdorff derivative.  
For $x \in X_n$ and $u,v \in \polys$,
\begin{equation}
\langle x u, v \rangle = \langle u, \twist_x v \rangle.
\end{equation}
In particular, for $P, Q \in \polys$,
$$\langle P, Q \rangle = P( \twist_{X_n} )\tau( Q ) \coeff_{x_1 = x_2 = \cdots = x_n =0}$$
where $\tau$ be the operator on $\polys$ that any monomial is sent to
the monomial obtained by reading its entries  from right-to-left.

%
%


\begin{defn}
Let $X_n=\{x_1, x_2, \cdots, x_n\}$ be a finite noncommuting alphabet. The 
harmonics of the symmetric group, with respect to the twisted derivative, are
defined as the space of solutions for the system of PDEs:
$$f(\twist_{X_n}) Q(X_n) =0$$ for all $f \in \NCSI_n$ without constant terms. We denote them 
as $\NCHar$.
\end{defn}

\begin{lemma} 
\[
\NCHar = \NCSIdeal^{\perp}~.
\]
\end{lemma}
\begin{proof} 
By definition, it is immediate that $\NCHar \subseteq \NCSIdeal^{\perp}.$

Suppose that $f \in  \NCSIdeal^{\perp}$. Then, for all $P$ in $\NCSIdeal$,
$$P(\twist_{X_n}) f(X_n) |_{x_1 = x_2 = \cdots = x_n =0} =0.$$
 We claim that this implies that $P(\twist_{X_n}) f(X_n) =0.$ Suppose
this is not the case. Then, let $u$ be the smallest monomial in
lexicographic order  that appears with nonzero coefficient in $P(\twist_{X_n}) f(X_n) $. 
Since $P$ is in  the left ideal $\NCSIdeal$, so is $u P$. But
by construction $(u P)(\twist_{X_n}) \, f(X_n) \ne 0.$
Contradiction. Hence, $\NCHar = \NCSIdeal^{\perp}.$
\end{proof}

Now we proceed as we did before in the case of $\MH_n$.  We have shown that
$\polys = \NCHar \oplus \NCSIdeal$ and by Proposition \ref{prop:basisIdl}
any $G(X_n) \in \polys$ can be expressed uniquely as
$$G(X_n) = f(X_n) + \sum_A P_A(X_n) \M_A[X_n]$$ 
where the sum is over nonsplitable set partitions $A$, $f(X_n) \in \NCHar$, and
$P_A(X_n) \in \polys$ is of degree strictly smaller than the degree of $G$.  This
procedure can be repeated recursively on $P_A(X_n)$ and the products of $\M_A[X_n]$
expanded in terms of other basis elements for $\NCSI$ so that
$$G(X_n) = \sum_{A} f_A(X_n) \M_A[X_n]$$
where the sum is over all set partitions $A$ of size smaller than or equal to the degree
of $G(X_n)$ and
each $f_A(X_n) \in \NCHar$.  This reduction is unique
and so the map $\opsi : \NCHar \otimes \NCSI_n \rightarrow \polys$ defined as
the linear extension of the map
$\opsi(f(X_n) \otimes P[X_n]) = f(X_n) P[X_n]$ is an isomorphism
of vector spaces.

\begin{prop}
As graded $\s_n$-modules,
\[
\NCHar \simeq \polys/\NCSIdeal
\]
\end{prop}
\begin{proof}
For each $G(X_n) \in \polys$, we know by the previous discussion that
since $\polys  = \NCHar \oplus \NCSIdeal$, there is a unique expression
$G(X_n) = f(X_n) +\sum_A P_A(X_n) \M_A[X_n]$ ($A$ are nonsplitable) 
and hence we have
that $G(X_n) \equiv f(X_n)~~(mod~\NCSIdeal)$.  Since for each non-empty
set partition $A$, $\sigma( P_A(X_n) \M_A[X_n] ) =
\sigma(P_A(X_n)) \M_A[X_n] \in \NCSIdeal$ we have that
$\sigma G(X_n) \equiv \sigma f(X_n)~~(mod~\NCSIdeal)$.
\end{proof}

%

\begin{theorem} \label{NCChevalley}
As graded $\s_n$-modules,
\[
\NCHar \otimes \NCSI_n  \simeq \polys~.
\]
\end{theorem}

\begin{proof}
Now consider the $\s_n$ action on the space $\NCHar \otimes \NCSI_n$.  Clearly
$\NCSI_n$ is an $\s_n$ module since each element is $\s_n$ invariant.  
For $f(X_n) \in \NCHar$ and for each $P(X_n) \in \NCSI_n$, $P(\twist_{X_n}) f(X_n) = 0$.
Since $\sigma( P(\twist_{X_n}) f(X_n) ) = \sigma(P(\twist_{X_n})) \sigma( f(X_n) )
= P(\twist_{X_n}) \sigma(f(X_n))$, hence $\sigma f(X_n) \in \NCHar$ also.

Now for $G(X_n) = \opsi \left( \sum_{A} f_A(X_n) \otimes \M_A[X_n] \right)$, we notice
that
\begin{equation}
\sigma G(X_n) = \sigma\left( \sum_{A} f_A(X_n) \M_A[X_n] \right) = 
\sum_{A} \sigma(f_A(X_n)) \M_A[X_n]
\end{equation}
and hence $\sigma \circ \opsi = \opsi \circ (\sigma \otimes \sigma)$ and therefore
the isomorphism holds on the level of $\s_n$-modules.
\end{proof}




We provide below a table of values of the dimensions of
the graded component of degree $k$ of the
space $\polys / \NCSIdeal$.
This is computed using the table of values in section 
\ref{sec:Wolfsummary} and the formula given in
Corollary \ref{cor:dimquo}.  The rows in the table below correspond to
the coefficients of $C_n(q)$.

\medskip

$$
\begin {array}{cccccccccc} 
n\slash k&\vline&0&1&2&3&4&5&6&7\\
\hline
1&\vline&1&0&0&0&0&0&0&0\\
2&\vline&1&1&1&1&1&1&1&1\\
3&\vline&1&2&5&13&34&89&233&610\\
4&\vline&1&3&11&42&162&627&2430&9423\\
5&\vline&1&4&19&93&459&2273&11274&55964\\
6&\vline&1&5&29&172&1026&6134&36712&219847\\
7&\vline&1&6&41&285&1989&13901&97215&680079\\
8&\vline&1&7&55&438&3498&27962&223604&1788406\end {array}
$$

\medskip

\end{section}

 \end{document}